\documentclass[hidelinks,onefignum,onetabnum]{siamart250211}

\usepackage{amsfonts}
\usepackage{graphicx}
\usepackage{epstopdf}
\ifpdf
  \DeclareGraphicsExtensions{.eps,.pdf,.png,.jpg}
\else
  \DeclareGraphicsExtensions{.eps}
\fi

\usepackage{wrapfig} 
\usepackage{ulem}
\usepackage{multirow}
\usepackage{graphicx}
\usepackage{amsmath}
\usepackage{amssymb}
\usepackage{caption}
\usepackage{subcaption}
\usepackage{tikz}
\usetikzlibrary{calc}
\usepackage{tabularx, makecell, booktabs}

\usepackage{mystylefile}
\usepackage{float}
\newsiamremark{remark}{Remark}
\newsiamremark{hypothesis}{Hypothesis}
\crefname{hypothesis}{Hypothesis}{Hypotheses}
\newsiamthm{claim}{Claim}
\newsiamremark{fact}{Fact}
\crefname{fact}{Fact}{Facts}

\headers{Randomized Block Low-Rank Matrix Compression by Tagging}{K.J. Pearce, A. Yesypenko, J. Levitt, P.-G. Martinsson}

\title{Randomized Block Low-Rank Matrix Compression\\ by Tagging\thanks{Submitted to the editors May 7, 2025.
\funding{This work was funded by the National Science Foundation (DMS-2401889, DMS-1952735, DMS-2313434), Office of Naval Research (N00014-18-1-2354), and the Department of Energy (DE-SC0022251).}}}

\author{Katherine J. Pearce\thanks{Department of Mathematics \& Oden Institute, University of Texas at Austin
  (\email{katherine.pearce@austin.utexas.edu}, \email{pgm@oden.utexas.edu}).}
\and Anna Yesypenko\thanks{Oden Institute, University of Texas at Austin
  (\email{annayesy@utexas.edu}, \email{jlevitt@utexas.edu}).}
\and James Levitt\footnotemark[3] \and Per-Gunnar Martinsson\footnotemark[2]}

\usepackage{amsopn}

\begin{document}

\maketitle

\begin{abstract}
In this work, we present randomized compression algorithms for flat rank-structured matrices with shared bases, termed uniform Block Low-Rank (BLR) matrices. Our main contribution is a technique called tagging, which improves upon the efficiency of existing algorithms for basis matrix computation while preserving accuracy. Tagging operates on the matrix using matrix-vector products of the matrix and its adjoint, making it suitable for black-box environments where accessing individual matrix entries is computationally expensive or infeasible.
 We show tagging requires a constant number of matrix-vector products coupled with linear post-processing; crucially, the asymptotic pre-factors in tagging depend only on the rank parameter and the underlying problem geometry.
We also establish a theoretical connection between the optimal construction of tagging matrices and projective varieties in algebraic geometry, suggesting a hybrid numeric-symbolic avenue of future work.
To validate our approach, we apply tagging to compress uniform BLR matrices arising from the discretization of integral and partial differential equations. Empirical results show that tagging outperforms alternative compression techniques, significantly reducing both the number of required matrix-vector products and overall computational time. These findings highlight the practicality and scalability of tagging as an efficient method for flat rank-structured matrices in scientific computing.
\end{abstract}



\section{Introduction}
\label{sec:intro}

In scientific computing and data science, many applications involve matrices that are dense but ``data-sparse,'' admitting certain low-rank approximations that compress the matrices while preserving critical information.
Algorithms to compress data-sparse matrices can achieve better performance by invoking \textit{rank structure}, where the input matrices are tessellated into blocks that are either small enough in size to apply dense algorithms or are of low numerical rank.

\begin{figure}
\centering
\begin{subfigure}{.5\textwidth}
  \centering
  \includegraphics[width=.35\linewidth]{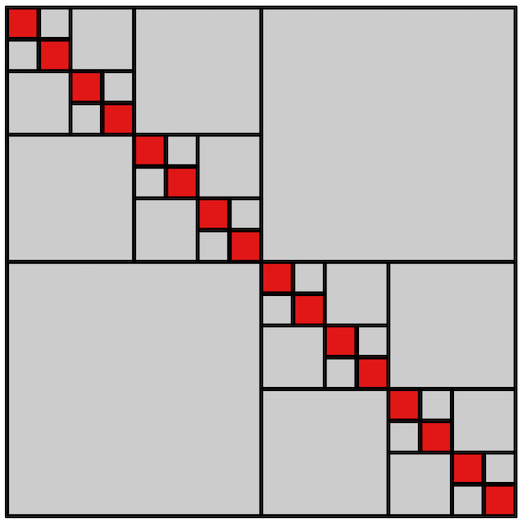}
  \caption{Hierarchical matrix format}
  \label{fig:sub1}
\end{subfigure}%
\begin{subfigure}{.5\textwidth}
  \centering
  \includegraphics[width=.35\linewidth]{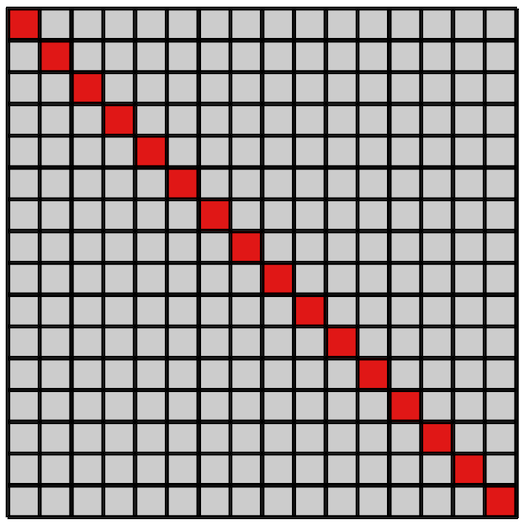}
  \caption{Flat matrix format}
  \label{fig:sub2}
\end{subfigure}
\caption{An illustration of hierarchical vs. flat rank-structured matrix formats. In (a), the rank-structure format uses a hierarchy of tessellations of the underlying problem domain. In (b), the flat format corresponds to a single tessellation of the domain; \textcolor{black}{more specifically, the flat rank-structured matrix in (b) corresponds to the leaf level of the tessellation of the hierarchical matrix in (a). As such, flat formats are the key building blocks of hierarchical formats and an important area of investigation in their own right, cf. \cite{amestoy2019}.}}
\label{fig:intro_fiat_vs_hierarchical}
\end{figure}

Many rank-structured matrix formats have been successfully utilized in engineering and data science applications.
These formats are typically classified as either hierarchical \cite{bebendorf2008hierarchical, 2010_borm_book, chandrasekaran2006fast, gillman2012direct, gorman2019robust, hackbusch1999sparse} or flat \cite{amestoy2015, amestoy2017, ashcraft2021, HighamMary2021, Jeannerod2019, Mary2017}; illustrative examples of hierarchical vs. flat rank-structure formats are shown in Fig.~\ref{fig:intro_fiat_vs_hierarchical}.
Both formats can be further characterized by either a weak or strong admissibility criterion.
In weakly admissible formats (as in Fig.~\ref{fig:intro_fiat_vs_hierarchical}), every off-diagonal block is admissible, or treated as low-rank, whereas admissible blocks of strongly admissible formats (cf.~Fig.~\ref{fig:blr2}) correspond to ``far-field" interactions (as in, e.g., the Fast Multipole Method \cite{greengard1987fast}). 

The last defining feature of rank-structured matrices most typically pertains to hierarchical formats, namely whether all admissible blocks in the same block-row or block-column are well-approximated by the same low-rank basis matrices. 
Frequently, basis matrix computations can be accelerated without impacting the accuracy of the low-rank approximation by invoking a shared basis assumption, namely that the row or column spaces of all admissible blocks within the same block-row or block-column are spanned by the same basis. 
A shared basis assumption is utilized in the strongly admissible uniform $\mathcal{H}^1$-matrix \cite{hackbusch1999sparse, lin2011fast} and $\mathcal{H}^2$-matrix \cite{giebermann2001multilevel, Hackbusch2000} formats, while the HSS \cite{chandrasekaran2006fast, xia2010superfast} and HBS \cite{gillman2012direct,levitt2022linear} formats use shared bases under weak admissibility. 


\textcolor{black}{In this manuscript, we focus on a flat format termed block low-rank (BLR) \cite{amestoy2015, amestoy2017, amestoy2019, HighamMary2021}.
BLR matrices have been successfully utilized in a variety of applications, including sparse direct solvers \cite{Amestoy2019_performance, Pichon2018, pichon2018_supernodal}, mathematical modeling \cite{Abdulah2018, akbudak2017, Cao2020, Shantsev2017}, multifrontal solvers \cite{Weisbecker2013, Mary2017}, boundary integral equations \cite{alHarthi2020, Sergent2016}, and deep neural networks \cite{chen2022pixelated, dao2022monarch, lee2024blast, lee2024differentiable, vaswani2017attention}.
Although flat matrix formats like BLR cannot be compressed in linear time like hierarchical formats, flat formats still enjoy many practical benefits that should be considered when choosing an appropriate rank structure format to represent a given problem.
For example, algorithms to construct or compress hierarchical matrices are difficult to parallelize (cf. \cite{Shantsev2017}), whereas the performance of algorithms for BLR matrices has been optimized to reap the full benefits of parallelization and GPU acceleration \cite{Amestoy2024, Amestoy2019_performance, beaumont:hal-02900244, Cao2020, Charara2018, Ida2018, Mary2017, Mondal2023, Pei2019, pichon2018_supernodal, Sergent2016}.
Additionally, much of the recent work on the BLR format has focused on proving theoretical guarantees, including stability analysis of solving BLR linear systems \cite{HighamMary2021, Mary2017}, complexity analysis of BLR algorithms \cite{amestoy2017, Amestoy2019_performance,  Mary2017}, and fast matrix arithmetic with BLR matrices \cite{amestoy2023mixed, apriansyah2022parallel, ashcraft2021, Ida2022, Jeannerod2019, Ohshima2022}.
In our work, we consider BLR matrices under a shared basis assumption as recently introduced in \cite{ashcraft2021}, which we refer to as uniform BLR matrices throughout this work.
Uniform BLR matrices have significantly smaller storage costs than BLR matrices, and the computational complexity of downstream tasks like LU factorization or inversion is greatly reduced with the shared basis assumption \cite{ashcraft2021, YenSDH23}.}

Though the BLR format has enjoyed much practical success and performance optimization, its utility in many applications has not yet been explored, particularly those in which the matrix entries cannot be directly accessed.
Rather, in these black-box problems, {or ``matrix-free'' environments},  the input matrix is assumed to be accessible only through a fast method to evaluate matrix-vector products, rather than through direct access to its entries.
Here, the goal is often matrix ``reconstruction'' in terms of low-rank basis matrices, enabling downstream operations (e.g., inversion or LU factorization) and simplifying formation of matrix products. 
These compressed matrix representations have broad applicability in scientific computing, for instance in deriving rank-structured representations of integral operators \cite{2005_martinsson_fastdirect, Yesypenko23_SkelFMM} or accelerating sparse direct solvers \cite{martinsson2019fast, yesypenko2023dissertation}. 

Randomized algorithms have been shown to be very effective in black-box problem environments, particularly \textit{randomized sketching}, in which the row and column spaces of an input matrix are approximated by analyzing how the matrix and its transpose act on tall thin matrices drawn from random matrix distributions \cite{halko2011finding, levitt2022linear, lin2011fast, 2008_martinsson_randomhudson}. 
Namely, suppose that $\ma \in \mathbb{R}^{N \times N}$ has some given rank structure where the numerical ranks of admissible blocks are upper bounded by $k \ll N$, but that $\ma$ is only accessible through some fast black-box algorithm.
In other words, given tall thin $\momega, \mpsi \in \mathbb{R}^{N \times r}$, $r=\mathcal{O}(k)$, we can quickly evaluate $\mY = \ma \momega$ and $\mZ = \ma^* \mpsi$.
Our goal is then to reconstruct $\ma$ 
as efficiently as possible, using only the information in the set $\{\mY,\momega,\mZ, \mpsi \}$
\cite{lin2011fast, martinsson2011fast, martinsson2016compressing}.

  Levitt and Martinsson \cite{levitt2022linear} introduced the first fully black-box linear-complexity randomized algorithm to compress (weakly admissible) HBS matrices, inspired by the linear-complexity HBS compression algorithm of \cite{martinsson2011fast}. Recent work \cite{Yesypenko23_RSRS} proposed a fully black-box algorithm for the simultaneous compression and LU factorization of general $\mathcal{H}^2$ matrices under strong admissibility, also achieving linear complexity. The present work focuses on the efficient black-box compression of uniform BLR matrices. Although flat formats have higher asymptotic complexity compared to hierarchical formats, they often feature modest constant factors in storage and can yield significant practical speed-ups.
  
We first adapt the block nullification algorithm of \cite{levitt2022linear,Yesypenko23_RSRS} for strongly admissible uniform BLR matrices. 
We then introduce a new randomized compression algorithm for uniform BLR matrices based on a method we refer to as \textit{tagging}, in which null spaces of small submatrices are computed to exclude contributions from inadmissible blocks in the random sketches $\mY$ and $\mZ$. 
As in \cite{levitt2022linear,Yesypenko23_RSRS}, the strategic exclusion of inadmissible blocks permits the same sample matrices $\mY$ and $\mZ$ to be used to compress every admissible block-row and block-column of the input matrix, which can be straightforwardly parallelized for optimized performance.
Additionally, our tagging method has a significantly smaller asymptotic pre-factor, and more critically, its sampling complexity for basis matrix computation is constant---independent of the problem size even in flat formats. 
We show that tagging improves on the performance of block nullification while maintaining approximation accuracy.

\textbf{Contributions.}
We propose two new  randomized compression schemes for flat rank-structured matrices with shared bases, termed uniform BLR.
The first is based on our extension of the block nullification method in \cite{levitt2022linear,Yesypenko23_RSRS} to uniform BLR matrices under a strong admissibility criterion.
The second is our main contribution, the novel method of tagging to compress uniform BLR matrices.
We draw a connection between the tagging matrices utilized in our method and Pl\"ucker coordinates in projective space, before presenting a highly efficient method of generating tagging matrices that works well in practice. 
We analyze the asymptotic complexities of both schemes and compare their performance in compressing matrices arising from discretizations of integral and partial differential equations, demonstrating the superior computational efficiency of tagging.

\textbf{Outline.} Our manuscript is structured as follows. Section~\ref{sec:bg} covers necessary linear algebra preliminaries and details the uniform BLR matrix format.
Section~\ref{sec:blocknull} illustrates our modification of the block nullification for basis matrix computations and analyzes its asymptotic complexity, and Section~\ref{sec:tagging} introduces the new method of tagging with detailed complexity analysis. Section~\ref{sec:tagmat_selection} outlines the relationship between tagging and projective varieties and provides practical strategies that yield high-quality random sketches with tagging.
Section~\ref{sec:full_comp_algs} finishes the compression procedure for uniform BLR matrices before showcasing numerical results in Section~\ref{sec:num_exp}.

\section{Preliminaries}
\label{sec:bg}

We briefly summarize necessary background for our work, following the presentation of~\cite{martinsson2016compressing}. We then outline the uniform BLR matrix format as in \cite{ashcraft2021}.

\subsection{Notation}

A vector $\xx \in \R^{n}$ is measured by the Euclidean norm $\Vert \xx \Vert = \left( \sum_i \vert x_i \vert^2 \right)^\frac{1}{2}$, and
a matrix $\ma \in \R^{m \times n}$ is equipped with the corresponding operator norm
$\Vert \ma \Vert = \sup_{\Vert \xx \Vert = 1} \Vert \ma \xx \Vert$. 
We let $[m]$ denote the integers $1,2,\ldots,m$.
We adopt the notation of Golub and Van Loan~\cite{golub2013matrix} to reference submatrices;
namely, if $\ma$ is an $m \times n$ matrix,
and $I = \lbrack i_1, i_2, \dots, i_k \rbrack \subset [m]$
and $J = \lbrack j_1, j_2, \dots, j_l \rbrack \subset [n]$ are (row and column, resp.) index sets,
then $\ma(I, J)$ denotes the $k \times l$ matrix
\begin{equation*}
    \ma(I, J)
    =
    \begin{bmatrix}
        \ma(i_1, j_1) & \ma(i_1, j_2) & \dots & \ma(i_1, j_l) \\
        \ma(i_2, j_1) & \ma(i_2, j_2) & \dots & \ma(i_2, j_l) \\
        \vdots        & \vdots        &       & \vdots        \\
        \ma(i_k, j_1) & \ma(i_k, j_2) & \dots & \ma(i_k, j_l) \\
    \end{bmatrix}.
\end{equation*}
The abbreviation $\ma(I,:)$ designates the submatrix
$\ma(I, [n])$,
and $\ma(:, J)$ is defined analogously.
The transpose of $\ma$ is given by $\ma^*$,
and a matrix $\mU$ is said to be \emph{orthonormal}
if its columns are orthonormal, $\mU^* \mU = \mI$.

\subsection{The QR factorization}
Every $m \times n$ matrix $\mB$ has a (full) QR factorization
\begin{equation}
\label{eq:QR}
\begin{array}{cccc}
    \mB        &=& \mQ        & \mR, \\
    m \times n & & m \times m & m \times n
\end{array}
\end{equation}
where
$\mQ$ is orthonormal
and
$\mR$ is upper-triangular.
If $\mB$ has rank $k$,
then it has a rank-$k$ partial QR factorization
given by
\begin{equation*}
\begin{array}{cccc}
    \mB_k      &=& \mQ_k       & \mR_k, \\
    m \times n & & m \times k  & k \times n
\end{array}
\end{equation*}
where $\mQ_k$ is orthonormal
and $\mR_k$ is upper-triangular.

\subsection{Functions for orthonormal bases}
\label{sec:preliminaries_on_basis}

For matrix $\mB$ of rank $\leq k$, the function outputting $\mQ$
with $k$ orthonormal columns spanning the column space of $\mB$ is
\begin{equation*}
\mQ = \col(\mB, k).
\end{equation*}
Analogously, for matrix $\mB$ with null space of dimension $\geq k$, we denote a function that returns a matrix $\mZ$ with $k$ orthonormal columns in the null space of $\mB$ by
\begin{equation*}
\mZ = \nullsp(\mB, k).
\end{equation*}
To implement $\col$ and $\nullsp$, a rank-revealing factorization, e.g., the SVD or column-pivoted QR, is usually necessary. However, for the particular input matrices $\mtx{B}$ arising in the present context, unpivoted QR can safely be used. (To be precise, the input matrices will be ``sample matrices'' in a randomized embedding, cf.~Section~\ref{sec:preliminaries_randomized}.)


\subsection{Randomized range-finding}
\label{sec:preliminaries_randomized}

Suppose $m \times n$ matrix $\mB$ can be accurately approximated by a rank-$k$ matrix, and we seek a matrix whose columns form an approximate orthonormal basis (ON-basis) for the column space of $\mB$.
Often referred to as range-finding, we seek an orthonormal matrix $\mQ$ such that
$\Vert \mB - \mQ \mQ^* \mB \Vert$ is small.
This problem can be addressed efficiently by \textit{randomized sketching}, where the column space of $\mB$ is approximated by analyzing how $\mB$ acts on matrices drawn from random distributions \cite{halko2011finding, Johnson1984ExtensionsOL, tropp2017practical}.
In general, randomized range-finding proceeds as follows:
\begin{enumerate}
    \item Choose a small integer $p$ for
        how much ``oversampling'' is done (e.g. $p=10$).
    \item Draw an $n \times (k+p)$ random matrix $\mG$ (e.g. Gaussian \cite{halko2011finding, Johnson1984ExtensionsOL}).
    \item Form the $m \times (k+p)$ random sketch $\mY = \mB \mG$.
    \item Compute $\mQ = \texttt{col}(\mY,k)$.
\end{enumerate}
Each column of $\mY$ is then a random linear combination of the columns of $\mB$, and the probability of obtaining an accurate column space approximation of $\mB$ with the column space of $\mY$ approaches 1 rapidly as $p$ increases;
notably, this probability depends only on $p$
(not on $m$ or $n$, or any other properties of $\mB$);
cf.~\cite{halko2011finding} and \cite[Section 11]{martinsson2020randomized}.
\textcolor{black}{In practice, when a target rank $k$ is unknown \textit{a priori}, this procedure can be implemented for adaptive rank determination using a prescribed error tolerance}.

\subsection{Block Low-Rank (BLR) Matrix Compression}
\label{sec:blr_main}

\input{fig_BLR2_matrix}

We focus on the randomized compression of $N \times N$ matrices $\ma$ admitting a \textit{block low-rank (BLR)} format.
BLR matrices are tessellated into $b$ row and $b$ column blocks according to a flat (vs. hierarchical) rank structure.
An example of a BLR matrix is illustrated in Fig.~\ref{fig:blr2}. 

Usually, strong admissibility is assumed for BLR matrices, as opposed to weak admissibility where every off-diagonal block is treated as low-rank (Fig.~\ref{fig:intro_fiat_vs_hierarchical}(b)), since the asymptotic complexity is the same as weak admissibility in flat formats \cite{ashcraft2021}.
Under a strong admissibility condition, the matrix blocks that correspond to the neighbors of a given box are also treated densely, i.e. not as low-rank (see Remark~\ref{remark:blradm}).
To determine the neighbors of a given block, the problem geometry is subdivided until each box (corresponding to a matrix block) reaches the desired cardinality \cite{George1973, Gilbert1986, ho2012fast, minden2017recursive}. 
The neighbors of a given matrix block are the blocks whose corresponding boxes share a boundary with the given box, including itself, for a total of up to $3^d$ neighbors in a $d$-dimensional geometry.
For example, in Fig.~\ref{fig:blr2}, the neighbor set of box 3 (outlined in magenta in the quasi-1D geometry) is given by $\mathcal{N}_3 = \{2,3,4\}$ (in blue).
The reader is referred to \cite{bebendorf2008hierarchical, hackbusch2015, 2005_martinsson_fastdirect, Mary2017} for more detailed treatments of rank-structured matrix admissibility, \textcolor{black}{including geometry-oblivious notions of admissibility as in \cite{yu2017geometry}}.


Analogous to the uniform $\mathcal{H}^1$ rank-structured matrix format \cite{hackbusch1999sparse, levitt2022randomized, lin2011fast}, we define a uniform BLR matrix $\ma$ as a BLR matrix for which low-rank blocks within the same block-row or block-column share the same bases of their row or column spaces.
Fig.~\ref{fig:unifbasis} illustrates this with $\ma$ from Fig.~\ref{fig:blr2}. 
To compress (strongly admissible) uniform BLR matrices, we compute basis matrices $\mU_{\ell}$ and $\mV_{m}$ for  low-rank block $\ma_{\ell,m}$ such that 
\begin{align}
\label{eq:uniformBLRblock}
    \underbrace{\ma_{\ell,m}}_{m \times m} = \underbrace{\mU_{\ell}}_{m \times k} \underbrace{\widetilde{\ma}_{\ell,m}}_{k \times k} \underbrace{\mV^*_{m}}_{k \times m},
\end{align}
where
\begin{align}
\begin{split}
\label{eq:uniformBLRbases}
    \mU_{\ell} &= \texttt{col} \left ( \ma(I_{\ell}, I_{i \in \mathcal{F}_{\ell}}), k \right ), \\
    \mV_{m} &= \texttt{col} \left ( \ma^*(I_{i \in \mathcal{F}_{m}}, I_{m}), k \right ).
\end{split}
\end{align}
Here, if $\mathcal{N}_{\ell}$ denotes the set of neighbors of box $\ell$, let $\mathcal{F}_{\ell} = [b] \setminus \mathcal{N}_{\ell} $ denote the set complement of $\mathcal{N}_\ell$,
and similarly for $\mathcal{F}_{m}$.
The shared bases assumption used in (\ref{eq:uniformBLRbases}) reduces the number of basis matrices from $\bigo(b^2)$ to $\bigo(b)$ for a $b \times b$ flat tessellation. 

When (\ref{eq:uniformBLRblock}) holds, we obtain a block factorization of an $N \times N$ uniform BLR matrix $\ma$ with $b$ block-rows/columns of size $m \times m$ (so that $N=bm$), given by
\begin{align}
\label{eq:uniformBLRfull}
    \underset{bm \times bm}{\ma} = \underset{bm \times bk} {\mU}\ \underset{bk \times bk} {\widetilde{\ma}}\ \underset{bk \times bm} {\mV}^* + \underset{bm \times bm} {\mB},
\end{align}
where
\begin{align}
\begin{split}
    \mU &= \mbox{diag}(\mU_1, \dots, \mU_b), \\
    \mV &= \mbox{diag}(\mV_1, \dots, \mV_b),
\end{split}
\end{align}
\begin{align}
    \widetilde{\ma} = \mU^* \ma \mV,
\end{align}
and $\mB$ is a block-sparse matrix defined (for strongly admissible $\ma$) as 
\begin{align}
\label{eq:discrepancymatrix}
\mB_{i, j} = \begin{cases} \ma_{i, j} - \mU_{i} \widetilde{\ma}_{i,j} \mV_{j}^*, &  i \in \mathcal{N}_{j} \\ \mathbf{0},& \text{otherwise}.
\end{cases}
\end{align}
The matrix $\mB$ given by (\ref{eq:discrepancymatrix}) represents the ``remainder'' of the inadmissible blocks of $\ma$ after their components spanned by the basis matrices have been peeled off. 

In general, uniform BLR compression can be accomplished by the following steps\footnote{Steps~\ref{comp2} and \ref{comp3} are interchangeable depending on the compression algorithm; see Section~\ref{sec:full_comp_algs}.}:
\begin{enumerate}[label=\textnormal{(\Roman*)}]
    \item \label{comp1} Compute basis matrices $\mU$ and $\mV$.
    \item \label{comp2} Compute matrix $\tilde{\ma} = \mU^* \ma \mV$.
    \item \label{comp3} Compute discrepancy matrix $\mB = \ma - \mU \widetilde{\ma} \mV^*$.
\end{enumerate}
The primary focus of this work is step~\ref{comp1}: we develop and compare randomized algorithms for computing basis matrices of uniform BLR matrices under strong admissibility conditions; \textcolor{black}{the only difference in reconstruction via our proposed methods is the method of basis construction (step~(I)).}
\textcolor{black}{We present two algorithms to compute basis matrices: the method of block nullification in Section~\ref{sec:blocknull} and our improved method of tagging in Section~\ref{sec:tagging}.
Crucially, these algorithms are suitable for black-box problem environments where we lack direct access to the matrix entries and their asymptotic complexity depends on the choice of $m$ (or $b$) and $k$; we discuss this in detail at the beginning of Section~\ref{sec:num_exp} after introducing each algorithm. }

\textcolor{black}{
Finally, we note that if matrix entries are readily available, steps \ref{comp1}-\ref{comp3} can be accomplished efficiently by the method of \cite{martinsson2011fast}.
Namely, the inadmissible matrix blocks  are explicitly subtracted off from $\ma$, so that the block factorization in (\ref{eq:uniformBLRfull}) takes a slightly different form and the low-rank basis can be computed with $\mathcal{O}(k+p)$ matrix-vector products (mat-vecs).
To see this, let $\mD \in \R^{N \times N}$ be the block sparse matrix of inadmissible blocks (corresponding to near-neighbor interactions, the red blocks in Fig.~\ref{fig:blr2}), with zeroes elsewhere. 
We contrast this with the matrix $\mB$ defined in (\ref{eq:discrepancymatrix}), which is the block-sparse remainder after the components of the inadmissible blocks spanned by the low-rank bases are subtracted off. 
In this case, the block factorization has the form
\begin{align}
    \ma = \mU \widehat{\ma} \mV^* + \mD,
\end{align}
where $\widehat{\ma} = \mU^*(\ma - \mD) \mV$.
In other words, we can compute the low-rank bases with $\mathcal{O}(k+p)$ mat-vecs with $\ma - \mD$ when we can access individual matrix entries.
}

\textcolor{black}{
When matrix entries are available, interpolative bases of the row and column spaces of admissible blocks are often preferable, instead of the orthonormal bases $\mU_1,\ldots,\mU_b$ and $\mV_1,\ldots,\mV_d$ that we compute in our work as given by (\ref{eq:uniformBLRbases}). 
Interpolative decompositions use actual rows and columns of a given matrix, often called skeletons, as approximate bases for its row or column space, where skeleton indices are typically chosen as the first $k$ pivots in column-pivoted QR (applied to random sketches of admissible blocks or their transposes, cf. \cite[Chapter 18]{martinsson2019fast}).
In this case, the complexity of the full uniform BLR reconstruction algorithm is dominated by the cost of accessing matrix entries, which is necessary not only to form $\ma-\mD$, but also to form the dense skeleton-to-skeleton interaction matrices comprising $\widehat{\ma}$.
As such, in this work, we only consider the black-box setting where matrix entries are not available, recommending that the existing method of \cite{martinsson2011fast} be used otherwise.
Additional implementation details for interpolative functions to replace \texttt{col} can be found in \cite[Chapter 18]{martinsson2019fast}.
}

\begin{remark}
\label{remark:blradm}
    \textcolor{black}{There are a few subtleties to point out related to admissibility. In our work, when we refer to weak or strong admissibility, we are referring to which blocks in the matrix are treated as low-rank, usually corresponding to some underlying problem geometry. However, there is also a notion of uniform BLR admissibility, or which blocks should be included in a shared basis. For the ease of explanation, we assume that all low-rank blocks in the same block-row/column share a basis for their row or column spaces. This is not necessary in general, cf. \cite{ashcraft2021}. In fact, in practice as noted in \cite{ashcraft2021}, inadmissible blocks (i.e. those treated densely) can often be compressed themselves (albeit to a lesser degree) for reduced storage costs and better performance in downstream tasks. }
\end{remark}

\begin{figure}[ht]
\centering
\begin{minipage}{0.22\linewidth}
    \fontsize{8}{\baselineskip}
    \def\n{5.25cm}
    \begin{tikzpicture}[level distance=\n/10,
        level 1/.style={sibling distance=\n/2},
        level 2/.style={sibling distance=\n/4},
        level 3/.style={sibling distance=\n/8},
        scale=0.5,]
        \begin{scope}[local bounding box=scope1, line width=1pt]

            \draw (0, 0) rectangle (\n, \n);

            \draw[fill=\lrcolor] (\n*2/8, \n*7/8) rectangle (\n*3/8, \n*8/8) node[midway] {};
            \draw[fill=\lrcolor] (\n*3/8, \n*7/8) rectangle (\n*4/8, \n*8/8) node[midway] {};
            \draw[fill=\unifcolscolor] (\n*4/8, \n*7/8) rectangle (\n*5/8, \n*8/8) node[midway] {};
            \draw[fill=\lrcolor] (\n*5/8, \n*7/8) rectangle (\n*6/8, \n*8/8) node[midway] {};
            \draw[fill=\lrcolor] (\n*6/8, \n*7/8) rectangle (\n*7/8, \n*8/8) node[midway] {};
            \draw[fill=\lrcolor] (\n*7/8, \n*7/8) rectangle (\n*8/8, \n*8/8) node[midway] {};

            \draw[fill=\lrcolor] (\n*3/8, \n*6/8) rectangle (\n*4/8, \n*7/8) node[midway] {};
            \draw[fill=\unifcolscolor] (\n*4/8, \n*6/8) rectangle (\n*5/8, \n*7/8) node[midway] {};
            \draw[fill=\lrcolor] (\n*5/8, \n*6/8) rectangle (\n*6/8, \n*7/8) node[midway] {};
            \draw[fill=\lrcolor] (\n*6/8, \n*6/8) rectangle (\n*7/8, \n*7/8) node[midway] {};
            \draw[fill=\lrcolor] (\n*7/8, \n*6/8) rectangle (\n*8/8, \n*7/8) node[midway] {};

            \draw[fill=\unifrowscolor] (\n*0/8, \n*5/8) rectangle (\n*1/8, \n*6/8) node[midway] {};
            \draw[fill=white] (\n*4/8, \n*5/8) rectangle (\n*5/8, \n*6/8) node[midway] {};
            \draw[fill=\unifrowscolor] (\n*5/8, \n*5/8) rectangle (\n*6/8, \n*6/8) node[midway] {};
            \draw[fill=\unifrowscolor] (\n*6/8, \n*5/8) rectangle (\n*7/8, \n*6/8) node[midway] {};
            \draw[fill=\unifrowscolor] (\n*7/8, \n*5/8) rectangle (\n*8/8, \n*6/8) node[midway] {};

            \draw[fill=\lrcolor] (\n*0/8, \n*4/8) rectangle (\n*1/8, \n*5/8) node[midway] {};
            \draw[fill=\lrcolor] (\n*1/8, \n*4/8) rectangle (\n*2/8, \n*5/8) node[midway] {};
            \draw[fill=\lrcolor] (\n*5/8, \n*4/8) rectangle (\n*6/8, \n*5/8) node[midway] {};
            \draw[fill=\lrcolor] (\n*6/8, \n*4/8) rectangle (\n*7/8, \n*5/8) node[midway] {};
            \draw[fill=\lrcolor] (\n*7/8, \n*4/8) rectangle (\n*8/8, \n*5/8) node[midway] {};

            \draw[fill=\lrcolor] (\n*0/8, \n*3/8) rectangle (\n*1/8, \n*4/8) node[midway] {};
            \draw[fill=\lrcolor] (\n*1/8, \n*3/8) rectangle (\n*2/8, \n*4/8) node[midway] {};
            \draw[fill=\lrcolor] (\n*2/8, \n*3/8) rectangle (\n*3/8, \n*4/8) node[midway] {};
            \draw[fill=\lrcolor] (\n*6/8, \n*3/8) rectangle (\n*7/8, \n*4/8) node[midway] {};
            \draw[fill=\lrcolor] (\n*7/8, \n*3/8) rectangle (\n*8/8, \n*4/8) node[midway] {};

            \draw[fill=\lrcolor] (\n*0/8, \n*2/8) rectangle (\n*1/8, \n*3/8) node[midway] {};
            \draw[fill=\lrcolor] (\n*1/8, \n*2/8) rectangle (\n*2/8, \n*3/8) node[midway] {};
            \draw[fill=\lrcolor] (\n*2/8, \n*2/8) rectangle (\n*3/8, \n*3/8) node[midway] {};
            \draw[fill=\lrcolor] (\n*3/8, \n*2/8) rectangle (\n*4/8, \n*3/8) node[midway] {};
            \draw[fill=\lrcolor] (\n*7/8, \n*2/8) rectangle (\n*8/8, \n*3/8) node[midway] {};

            \draw[fill=\lrcolor] (\n*0/8, \n*1/8) rectangle (\n*1/8, \n*2/8) node[midway] {};
            \draw[fill=\lrcolor] (\n*1/8, \n*1/8) rectangle (\n*2/8, \n*2/8) node[midway] {};
            \draw[fill=\lrcolor] (\n*2/8, \n*1/8) rectangle (\n*3/8, \n*2/8) node[midway] {};
            \draw[fill=\lrcolor] (\n*3/8, \n*1/8) rectangle (\n*4/8, \n*2/8) node[midway] {};
            \draw[fill=\unifcolscolor] (\n*4/8, \n*1/8) rectangle (\n*5/8, \n*2/8) node[midway] {};

            \draw[fill=\lrcolor] (\n*0/8, \n*0/8) rectangle (\n*1/8, \n*1/8) node[midway] {};
            \draw[fill=\lrcolor] (\n*1/8, \n*0/8) rectangle (\n*2/8, \n*1/8) node[midway] {};
            \draw[fill=\lrcolor] (\n*2/8, \n*0/8) rectangle (\n*3/8, \n*1/8) node[midway] {};
            \draw[fill=\lrcolor] (\n*3/8, \n*0/8) rectangle (\n*4/8, \n*1/8) node[midway] {};
            \draw[fill=\unifcolscolor] (\n*4/8, \n*0/8) rectangle (\n*5/8, \n*1/8) node[midway] {};
            \draw[fill=\lrcolor] (\n*5/8, \n*0/8) rectangle (\n*6/8, \n*1/8) node[midway] {};

            \draw[fill=\frcolor] (\n*1/8, \n*7/8) rectangle (\n*2/8, \n*8/8) node[midway] {};
            \draw[fill=\frcolor] (\n*2/8, \n*6/8) rectangle (\n*3/8, \n*7/8) node[midway] {};
            \draw[fill=\frcolor] (\n*3/8, \n*5/8) rectangle (\n*4/8, \n*6/8) node[midway] {};
            \draw[fill=\frcolor] (\n*4/8, \n*4/8) rectangle (\n*5/8, \n*5/8) node[midway] {};
            \draw[fill=\frcolor] (\n*5/8, \n*3/8) rectangle (\n*6/8, \n*4/8) node[midway] {};
            \draw[fill=\frcolor] (\n*6/8, \n*2/8) rectangle (\n*7/8, \n*3/8) node[midway] {};
            \draw[fill=\frcolor] (\n*7/8, \n*1/8) rectangle (\n*8/8, \n*2/8) node[midway] {};

            \draw[fill=\frcolor] (\n*0/8, \n*6/8) rectangle (\n*1/8, \n*7/8) node[midway] {};
            \draw[fill=\frcolor] (\n*1/8, \n*5/8) rectangle (\n*2/8, \n*6/8) node[midway] {};
            \draw[fill=\frcolor] (\n*2/8, \n*4/8) rectangle (\n*3/8, \n*5/8) node[midway] {};
            \draw[fill=\frcolor] (\n*3/8, \n*3/8) rectangle (\n*4/8, \n*4/8) node[midway] {};
            \draw[fill=\frcolor] (\n*4/8, \n*2/8) rectangle (\n*5/8, \n*3/8) node[midway] {};
            \draw[fill=\frcolor] (\n*5/8, \n*1/8) rectangle (\n*6/8, \n*2/8) node[midway] {};
            \draw[fill=\frcolor] (\n*6/8, \n*0/8) rectangle (\n*7/8, \n*1/8) node[midway] {};

            \draw[fill=\frcolor] (\n*0/8, \n*7/8) rectangle (\n*1/8, \n*8/8) node[midway] {};
            \draw[fill=\frcolor] (\n*1/8, \n*6/8) rectangle (\n*2/8, \n*7/8) node[midway] {};
            \draw[fill=\frcolor] (\n*2/8, \n*5/8) rectangle (\n*3/8, \n*6/8) node[midway] {};
            \draw[fill=\frcolor] (\n*3/8, \n*4/8) rectangle (\n*4/8, \n*5/8) node[midway] {};
            \draw[fill=\frcolor] (\n*4/8, \n*3/8) rectangle (\n*5/8, \n*4/8) node[midway] {};
            \draw[fill=\frcolor] (\n*5/8, \n*2/8) rectangle (\n*6/8, \n*3/8) node[midway] {};
            \draw[fill=\frcolor] (\n*6/8, \n*1/8) rectangle (\n*7/8, \n*2/8) node[midway] {};
            \draw[fill=\frcolor] (\n*7/8, \n*0/8) rectangle (\n*8/8, \n*1/8) node[midway] {};
        \end{scope}
    \end{tikzpicture}
\end{minipage}
\raisebox{0mm}{
\begin{minipage}{0.7\linewidth}
\small{{
    \textit{For example, to compute $\ma_{3,5} = \mU_{3} \widetilde{\ma}_{3,5} \mV_{5}^*$ (white)}:
    \begin{align*}
        \mU_{3} &= \texttt{col} \left ( \begin{bmatrix} \ \mbox{\fcolorbox{black}{\unifrowscolor}{$\ma_{3,1}$}}  \ , \ \mbox{\fbox{$\ma_{3,5}$}} \ , \ \mbox{\fcolorbox{black}{\unifrowscolor}{$\ma_{3,6}$}} \ , \ \mbox{\fcolorbox{black}{\unifrowscolor}{$\ma_{3,7}$}} \ , \ \mbox{\fcolorbox{black}{\unifrowscolor}{$\ma_{3,8}$}} \textcolor{white}{j}  \end{bmatrix}, k \right ) \\
        \mV_{5} &= \texttt{col} \left ( \begin{bmatrix} \ \mbox{\fcolorbox{black}{\unifcolscolor}{$\ma_{1,5}$}} \ ; \  \mbox{\fcolorbox{black}{\unifcolscolor}{$\ma_{2,5}$}} \  ;  \ \mbox{\fbox{$\ma_{3,5}$}} \ ; \ \mbox{\fcolorbox{black}{\unifcolscolor}{$\ma_{7,5}$}} \  ; \   \mbox{\fcolorbox{black}{\unifcolscolor}{$\ma_{8,5}$}} \textcolor{white}{j} \end{bmatrix}^*, k \right ) \\
        \Rightarrow \widetilde{\ma}_{3,5} &= \mU_{3}^* \ma_{3,5} \mV_{5}.
    \end{align*}}}
\end{minipage}}
\caption{Basis matrices for a strongly-admissible uniform BLR matrix $\ma$ (Fig. \ref{fig:blr2})}
\label{fig:unifbasis}
\end{figure}

\section{Block Nullification in Uniform BLR Matrix Compression}
\label{sec:blocknull}

In this section, we modify the algorithm from \cite{levitt2022linear,Yesypenko23_RSRS} using ``block nullification'' to form random sketches of admissible matrix blocks.
These sketches are then used to compute basis matrices in an HBS representation according to the randomized rangefinding procedure of Section~\ref{sec:preliminaries_randomized}.
The algorithm is also fully black-box, accomplishing steps \ref{comp1}-\ref{comp3} without access to matrix entries.
Rather, it assumes access to fast matrix multiplication, so that sample matrices $\mY = \ma \momega$ and $\mZ = \ma^* \mpsi$ are formed efficiently given tall thin random test matrices $\momega, \mpsi \in \mathbb{R}^{N \times s}$ for $s = \mathcal{O}(k)$ for block-rank $k$.

The goal is then to ``reconstruct'' $\ma$ via steps~\ref{comp1}-\ref{comp3} using only $\mY, \mZ, \momega$, and $\mpsi$ (computed a priori).
However, each row of $\mY$, for instance, is a random linear combination of all columns of $\ma$ within a block-row, including the columns belonging to inadmissible blocks.
Block nullification yields ``clean'' random sketches from $\mY$ and $\mZ$ by excluding contributions from inadmissible blocks, without repeatedly applying $\ma$ or $\ma^*$ to tailored random test matrices to individually sketch admissible blocks. 

While block nullification is suitable for flat or hierarchical rank-structured formats, its performance has not been investigated for flat formats.
As such, we first modify block nullification for strongly admissible uniform BLR matrices.
We then discuss its asymptotic complexity to emphasize its larger pre-factor versus tagging. 

\subsection{Block nullification for strongly-admissible uniform BLR matrices}
\label{sec:blocknull_overview}

We begin with an illustrative example of the block nullification technique applied to the strongly-admissible uniform BLR matrix $\ma \in \mathbb{R}^{N \times N}$ from Fig.~\ref{fig:blr2}, tessellated into $b \times b$ blocks each of size $m \times m$ and target rank $k$.
Let $r = k+p$ for small $p$ (e.g. $p=10$), and let $\momega, \mpsi \in \mathbb{R}^{N \times s}$ be Gaussian test matrices with $s \geq r+3m$.

Suppose that we want to compute $\mU_{3} \in \mathbb{R}^{m \times k}$ as in Fig.~\ref{fig:unifbasis}, now using the randomized rangefinder procedure with a random sketch $\mY$ of the form 
\begin{center}
\resizebox{.75\linewidth}{!}{\begin{minipage}{\linewidth}
\begin{align*}
    \begin{array}{cccc}
      \mY   & = & \ma & \momega \\
       \begin{tikzpicture}[scale=0.79, line width=1pt]
            \def\n{4cm}
            \draw[fill=\lrcolor, dashed] (\n*0/8, \n*7/8) rectangle (\n*2/8, \n*8/8) node[midway] {};
            \draw[fill=\lrcolor, dashed] (\n*0/8, \n*6/8) rectangle (\n*2/8, \n*7/8) node[midway] {};
            \draw[fill=\lrcolor, dashed] (\n*0/8, \n*4/8) rectangle (\n*2/8, \n*5/8) node[midway] {};
            \draw[fill=\lrcolor, dashed] (\n*0/8, \n*3/8) rectangle (\n*2/8, \n*4/8) node[midway] {};
            \draw[fill=\lrcolor, dashed] (\n*0/8, \n*2/8) rectangle (\n*2/8, \n*3/8) node[midway] {};
            \draw[fill=\lrcolor, dashed] (\n*0/8, \n*1/8) rectangle (\n*2/8, \n*2/8) node[midway] {};
            \draw[fill=\lrcolor, dashed] (\n*0/8, \n*0/8) rectangle (\n*2/8, \n*1/8) node[midway] {};
            \draw[fill=\lrcolor] (\n*0/8, \n*5/8) rectangle (\n*2/8, \n*6/8) node[midway] {$\mY_{3}$};
\end{tikzpicture}  &  & \fontsize{8}{\baselineskip}
    \def\n{8cm}
    \begin{tikzpicture}[level distance=\n/10,
        level 1/.style={sibling distance=\n/2},
        level 2/.style={sibling distance=\n/4},
        level 3/.style={sibling distance=\n/8},
        scale=0.4,]
        \begin{scope}[local bounding box=scope1, line width=1pt]


            \draw[fill=\lrcolor, dashed] (\n*2/8, \n*7/8) rectangle (\n*3/8, \n*8/8) node[midway] {};
            \draw[fill=\lrcolor, dashed] (\n*3/8, \n*7/8) rectangle (\n*4/8, \n*8/8) node[midway] {};
            \draw[fill=\lrcolor, dashed] (\n*4/8, \n*7/8) rectangle (\n*5/8, \n*8/8) node[midway] {};
            \draw[fill=\lrcolor, dashed] (\n*5/8, \n*7/8) rectangle (\n*6/8, \n*8/8) node[midway] {};
            \draw[fill=\lrcolor, dashed] (\n*6/8, \n*7/8) rectangle (\n*7/8, \n*8/8) node[midway] {};
            \draw[fill=\lrcolor, dashed] (\n*7/8, \n*7/8) rectangle (\n*8/8, \n*8/8) node[midway] {};

            \draw[fill=\lrcolor, dashed] (\n*3/8, \n*6/8) rectangle (\n*4/8, \n*7/8) node[midway] {};
            \draw[fill=\lrcolor, dashed] (\n*4/8, \n*6/8) rectangle (\n*5/8, \n*7/8) node[midway] {};
            \draw[fill=\lrcolor, dashed] (\n*5/8, \n*6/8) rectangle (\n*6/8, \n*7/8) node[midway] {};
            \draw[fill=\lrcolor, dashed] (\n*6/8, \n*6/8) rectangle (\n*7/8, \n*7/8) node[midway] {};
            \draw[fill=\lrcolor, dashed] (\n*7/8, \n*6/8) rectangle (\n*8/8, \n*7/8) node[midway] {};

            \draw[fill=\unifrowscolor] (\n*0/8, \n*5/8) rectangle (\n*1/8, \n*6/8) node[midway] {};
            \draw[fill=\unifrowscolor] (\n*4/8, \n*5/8) rectangle (\n*5/8, \n*6/8) node[midway] {};
            \draw[fill=\unifrowscolor] (\n*5/8, \n*5/8) rectangle (\n*6/8, \n*6/8) node[midway] {};
            \draw[fill=\unifrowscolor] (\n*6/8, \n*5/8) rectangle (\n*7/8, \n*6/8) node[midway] {};
            \draw[fill=\unifrowscolor] (\n*7/8, \n*5/8) rectangle (\n*8/8, \n*6/8) node[midway] {};

            \draw[fill=\lrcolor, dashed] (\n*0/8, \n*4/8) rectangle (\n*1/8, \n*5/8) node[midway] {};
            \draw[fill=\lrcolor, dashed] (\n*1/8, \n*4/8) rectangle (\n*2/8, \n*5/8) node[midway] {};
            \draw[fill=\lrcolor, dashed] (\n*5/8, \n*4/8) rectangle (\n*6/8, \n*5/8) node[midway] {};
            \draw[fill=\lrcolor, dashed] (\n*6/8, \n*4/8) rectangle (\n*7/8, \n*5/8) node[midway] {};
            \draw[fill=\lrcolor, dashed] (\n*7/8, \n*4/8) rectangle (\n*8/8, \n*5/8) node[midway] {};

            \draw[fill=\lrcolor, dashed] (\n*0/8, \n*3/8) rectangle (\n*1/8, \n*4/8) node[midway] {};
            \draw[fill=\lrcolor, dashed] (\n*1/8, \n*3/8) rectangle (\n*2/8, \n*4/8) node[midway] {};
            \draw[fill=\lrcolor, dashed] (\n*2/8, \n*3/8) rectangle (\n*3/8, \n*4/8) node[midway] {};
            \draw[fill=\lrcolor, dashed] (\n*6/8, \n*3/8) rectangle (\n*7/8, \n*4/8) node[midway] {};
            \draw[fill=\lrcolor, dashed] (\n*7/8, \n*3/8) rectangle (\n*8/8, \n*4/8) node[midway] {};

            \draw[fill=\lrcolor, dashed] (\n*0/8, \n*2/8) rectangle (\n*1/8, \n*3/8) node[midway] {};
            \draw[fill=\lrcolor, dashed] (\n*1/8, \n*2/8) rectangle (\n*2/8, \n*3/8) node[midway] {};
            \draw[fill=\lrcolor, dashed] (\n*2/8, \n*2/8) rectangle (\n*3/8, \n*3/8) node[midway] {};
            \draw[fill=\lrcolor, dashed] (\n*3/8, \n*2/8) rectangle (\n*4/8, \n*3/8) node[midway] {};
            \draw[fill=\lrcolor, dashed] (\n*7/8, \n*2/8) rectangle (\n*8/8, \n*3/8) node[midway] {};

            \draw[fill=\lrcolor, dashed] (\n*0/8, \n*1/8) rectangle (\n*1/8, \n*2/8) node[midway] {};
            \draw[fill=\lrcolor, dashed] (\n*1/8, \n*1/8) rectangle (\n*2/8, \n*2/8) node[midway] {};
            \draw[fill=\lrcolor, dashed] (\n*2/8, \n*1/8) rectangle (\n*3/8, \n*2/8) node[midway] {};
            \draw[fill=\lrcolor, dashed] (\n*3/8, \n*1/8) rectangle (\n*4/8, \n*2/8) node[midway] {};
            \draw[fill=\lrcolor, dashed] (\n*4/8, \n*1/8) rectangle (\n*5/8, \n*2/8) node[midway] {};

            \draw[fill=\lrcolor, dashed] (\n*0/8, \n*0/8) rectangle (\n*1/8, \n*1/8) node[midway] {};
            \draw[fill=\lrcolor, dashed] (\n*1/8, \n*0/8) rectangle (\n*2/8, \n*1/8) node[midway] {};
            \draw[fill=\lrcolor, dashed] (\n*2/8, \n*0/8) rectangle (\n*3/8, \n*1/8) node[midway] {};
            \draw[fill=\lrcolor, dashed] (\n*3/8, \n*0/8) rectangle (\n*4/8, \n*1/8) node[midway] {};
            \draw[fill=\lrcolor, dashed] (\n*4/8, \n*0/8) rectangle (\n*5/8, \n*1/8) node[midway] {};
            \draw[fill=\lrcolor, dashed] (\n*5/8, \n*0/8) rectangle (\n*6/8, \n*1/8) node[midway] {};

            \draw[fill=\frcolor, dashed] (\n*1/8, \n*7/8) rectangle (\n*2/8, \n*8/8) node[midway] {};
            \draw[fill=\frcolor, dashed] (\n*2/8, \n*6/8) rectangle (\n*3/8, \n*7/8) node[midway] {};
            \draw[fill=\frcolor] (\n*3/8, \n*5/8) rectangle (\n*4/8, \n*6/8) node[midway] {};
            \draw[fill=\frcolor, dashed] (\n*4/8, \n*4/8) rectangle (\n*5/8, \n*5/8) node[midway] {};
            \draw[fill=\frcolor, dashed] (\n*5/8, \n*3/8) rectangle (\n*6/8, \n*4/8) node[midway] {};
            \draw[fill=\frcolor, dashed] (\n*6/8, \n*2/8) rectangle (\n*7/8, \n*3/8) node[midway] {};
            \draw[fill=\frcolor, dashed] (\n*7/8, \n*1/8) rectangle (\n*8/8, \n*2/8) node[midway] {};

            \draw[fill=\frcolor, dashed] (\n*0/8, \n*6/8) rectangle (\n*1/8, \n*7/8) node[midway] {};
            \draw[fill=\frcolor] (\n*1/8, \n*5/8) rectangle (\n*2/8, \n*6/8) node[midway] {};
            \draw[fill=\frcolor, dashed] (\n*2/8, \n*4/8) rectangle (\n*3/8, \n*5/8) node[midway] {};
            \draw[fill=\frcolor, dashed] (\n*3/8, \n*3/8) rectangle (\n*4/8, \n*4/8) node[midway] {};
            \draw[fill=\frcolor, dashed] (\n*4/8, \n*2/8) rectangle (\n*5/8, \n*3/8) node[midway] {};
            \draw[fill=\frcolor, dashed] (\n*5/8, \n*1/8) rectangle (\n*6/8, \n*2/8) node[midway] {};
            \draw[fill=\frcolor, dashed] (\n*6/8, \n*0/8) rectangle (\n*7/8, \n*1/8) node[midway] {};

            \draw[fill=\frcolor, dashed] (\n*0/8, \n*7/8) rectangle (\n*1/8, \n*8/8) node[midway] {};
            \draw[fill=\frcolor, dashed] (\n*1/8, \n*6/8) rectangle (\n*2/8, \n*7/8) node[midway] {};
            \draw[fill=\frcolor] (\n*2/8, \n*5/8) rectangle (\n*3/8, \n*6/8) node[midway] {};
            \draw[fill=\frcolor, dashed] (\n*3/8, \n*4/8) rectangle (\n*4/8, \n*5/8) node[midway] {};
            \draw[fill=\frcolor, dashed] (\n*4/8, \n*3/8) rectangle (\n*5/8, \n*4/8) node[midway] {};
            \draw[fill=\frcolor, dashed] (\n*5/8, \n*2/8) rectangle (\n*6/8, \n*3/8) node[midway] {};
            \draw[fill=\frcolor, dashed] (\n*6/8, \n*1/8) rectangle (\n*7/8, \n*2/8) node[midway] {};
            \draw[fill=\frcolor, dashed] (\n*7/8, \n*0/8) rectangle (\n*8/8, \n*1/8) node[midway] {};

            \draw[fill=\frcolor] (\n*2/8, \n*5/8) rectangle (\n*3/8, \n*6/8) node[midway] {};
            \draw[fill=\frcolor] (\n*1/8, \n*5/8) rectangle (\n*2/8, \n*6/8) node[midway] {};
            \draw[fill=\unifrowscolor] (\n*0/8, \n*5/8) rectangle (\n*1/8, \n*6/8) node[midway] {};
            \draw[fill=\frcolor] (\n*3/8, \n*5/8) rectangle (\n*4/8, \n*6/8) node[midway] {};
            \draw[fill=\unifrowscolor] (\n*4/8, \n*5/8) rectangle (\n*5/8, \n*6/8) node[midway] {};
            \draw[fill=\unifrowscolor] (\n*5/8, \n*5/8) rectangle (\n*6/8, \n*6/8) node[midway] {};
            \draw[fill=\unifrowscolor] (\n*6/8, \n*5/8) rectangle (\n*7/8, \n*6/8) node[midway] {};
            \draw[fill=\unifrowscolor] (\n*7/8, \n*5/8) rectangle (\n*8/8, \n*6/8) node[midway] {};
        \end{scope}
    \end{tikzpicture}
 &  \begin{tikzpicture}[scale=0.79, line width=1pt]
            \def\n{4cm}
            \draw[fill=\unifrowscolor] (\n*0/8, \n*7/8) rectangle (\n*2/8, \n*8/8) node[midway] {$\momega_1$};
            \draw[fill=\frcolor] (\n*0/8, \n*6/8) rectangle (\n*2/8, \n*7/8) node[midway] {$\momega_2$};
            \draw[fill=\frcolor] (\n*0/8, \n*5/8) rectangle (\n*2/8, \n*6/8) node[midway] {$\momega_{3}$};
            \draw[fill=\frcolor] (\n*0/8, \n*4/8) rectangle (\n*2/8, \n*5/8) node[midway] {$\momega_{4}$};
            \draw[fill=\unifrowscolor] (\n*0/8, \n*3/8) rectangle (\n*2/8, \n*4/8) node[midway] {$\momega_{5}$};
            \draw[fill=\unifrowscolor] (\n*0/8, \n*2/8) rectangle (\n*2/8, \n*3/8) node[midway] {$\momega_{6}$};
            \draw[fill=\unifrowscolor] (\n*0/8, \n*1/8) rectangle (\n*2/8, \n*2/8) node[midway] {$\momega_7$};
            \draw[fill=\unifrowscolor] (\n*0/8, \n*0/8) rectangle (\n*2/8, \n*1/8) node[midway] {$\momega_{8}$};
\end{tikzpicture}
    \end{array}
\end{align*} 
\end{minipage}}
\end{center}

\noindent where $\mY_{3} = \mY(I_{3},:)$, and $\momega_i = \momega(I_i,:)$ for $i = 1,\ldots,8$ are color-coded according to their respective block-factors in the block-row $\ma_{3} = \ma(I_{3},:)$.
As in Fig.~\ref{fig:unifbasis}, the blue blocks are the admissible blocks in $\ma_{3}$ whose columns will be approximately spanned by $\mU_{3}$ computed with $\mY_{3}$, whereas the red blocks of $\ma$ are inadmissible; our goal is to exclude their contributions from the randomized sample of $\ma_{3}$ held in $\mY_{3}$.

Since $\momega^{(3)}:=\momega([I_{2},I_{3},I_{4}],:) $ is of size $3m \times s$, it has a null space of dimension at least $s - 3m \geq r$.
We then compute a set of $r$ orthonormal vectors in its null space: 
\begin{align}
\label{eq:bn_nullsp}
    \underset{s \times r}{\mP^{(3)}} = \nullsp(\momega^{(3)}, r) := \nullsp(\begin{bmatrix}
        \momega_2 \\
        \momega_{3} \\
        \momega_{4}
    \end{bmatrix}, r) \ \ \Rightarrow \ \  \momega^{(3)} \mP^{(3)} = \begin{bmatrix}
        \momega_2 \\
        \momega_{3} \\
        \momega_{4}
    \end{bmatrix} \mP^{(3)} = \underset{3m \times r}{\boldsymbol{0}}.
\end{align}
Thus, we can obtain the desired sample of the admissible blocks of $\ma_{3}$ inexpensively from $\mY$ via  $\mY_{3} \mP^{(3)}$ (blue):
\begin{center}
\resizebox{.75\linewidth}{!}{\begin{minipage}{\linewidth}
\begin{align*}
        \begin{array}{cccc}
      \mY \mP^{(3)}   & = & \ma & \momega \mP^{(3)} \\
       \begin{tikzpicture}[scale=0.8, line width=1pt]
            \def\n{4cm}
            \draw[fill=\lrcolor, dashed] (\n*0/8, \n*7/8) rectangle (\n*0.66/8, \n*8/8) node[midway] {};
            \draw[fill=\lrcolor, dashed] (\n*0/8, \n*6/8) rectangle (\n*0.66/8, \n*7/8) node[midway] {};
            \draw[fill=\lrcolor, dashed] (\n*0/8, \n*4/8) rectangle (\n*0.66/8, \n*5/8) node[midway] {};
            \draw[fill=\lrcolor, dashed] (\n*0/8, \n*3/8) rectangle (\n*0.66/8, \n*4/8) node[midway] {};
            \draw[fill=\lrcolor, dashed] (\n*0/8, \n*2/8) rectangle (\n*0.66/8, \n*3/8) node[midway] {};
            \draw[fill=\lrcolor, dashed] (\n*0/8, \n*1/8) rectangle (\n*0.66/8, \n*2/8) node[midway] {};
            \draw[fill=\lrcolor, dashed] (\n*0/8, \n*0/8) rectangle (\n*0.66/8, \n*1/8) node[midway] {};
            \draw[fill=\unifrowscolor] (\n*0/8, \n*5/8) rectangle (\n*0.66/8, \n*6/8) node[midway] {};
\end{tikzpicture}  &  & \fontsize{8}{\baselineskip}
    \def\n{8cm}
    \begin{tikzpicture}[level distance=\n/10,
        level 1/.style={sibling distance=\n/2},
        level 2/.style={sibling distance=\n/4},
        level 3/.style={sibling distance=\n/8},
        scale=0.4,]
        \begin{scope}[local bounding box=scope1, line width=1pt]


            \draw[fill=\lrcolor, dashed] (\n*2/8, \n*7/8) rectangle (\n*3/8, \n*8/8) node[midway] {};
            \draw[fill=\lrcolor, dashed] (\n*3/8, \n*7/8) rectangle (\n*4/8, \n*8/8) node[midway] {};
            \draw[fill=\lrcolor, dashed] (\n*4/8, \n*7/8) rectangle (\n*5/8, \n*8/8) node[midway] {};
            \draw[fill=\lrcolor, dashed] (\n*5/8, \n*7/8) rectangle (\n*6/8, \n*8/8) node[midway] {};
            \draw[fill=\lrcolor, dashed] (\n*6/8, \n*7/8) rectangle (\n*7/8, \n*8/8) node[midway] {};
            \draw[fill=\lrcolor, dashed] (\n*7/8, \n*7/8) rectangle (\n*8/8, \n*8/8) node[midway] {};

            \draw[fill=\lrcolor, dashed] (\n*3/8, \n*6/8) rectangle (\n*4/8, \n*7/8) node[midway] {};
            \draw[fill=\lrcolor, dashed] (\n*4/8, \n*6/8) rectangle (\n*5/8, \n*7/8) node[midway] {};
            \draw[fill=\lrcolor, dashed] (\n*5/8, \n*6/8) rectangle (\n*6/8, \n*7/8) node[midway] {};
            \draw[fill=\lrcolor, dashed] (\n*6/8, \n*6/8) rectangle (\n*7/8, \n*7/8) node[midway] {};
            \draw[fill=\lrcolor, dashed] (\n*7/8, \n*6/8) rectangle (\n*8/8, \n*7/8) node[midway] {};

            \draw[fill=\unifrowscolor] (\n*0/8, \n*5/8) rectangle (\n*1/8, \n*6/8) node[midway] {};
            \draw[fill=\unifrowscolor] (\n*4/8, \n*5/8) rectangle (\n*5/8, \n*6/8) node[midway] {};
            \draw[fill=\unifrowscolor] (\n*5/8, \n*5/8) rectangle (\n*6/8, \n*6/8) node[midway] {};
            \draw[fill=\unifrowscolor] (\n*6/8, \n*5/8) rectangle (\n*7/8, \n*6/8) node[midway] {};
            \draw[fill=\unifrowscolor] (\n*7/8, \n*5/8) rectangle (\n*8/8, \n*6/8) node[midway] {};

            \draw[fill=\lrcolor, dashed] (\n*0/8, \n*4/8) rectangle (\n*1/8, \n*5/8) node[midway] {};
            \draw[fill=\lrcolor, dashed] (\n*1/8, \n*4/8) rectangle (\n*2/8, \n*5/8) node[midway] {};
            \draw[fill=\lrcolor, dashed] (\n*5/8, \n*4/8) rectangle (\n*6/8, \n*5/8) node[midway] {};
            \draw[fill=\lrcolor, dashed] (\n*6/8, \n*4/8) rectangle (\n*7/8, \n*5/8) node[midway] {};
            \draw[fill=\lrcolor, dashed] (\n*7/8, \n*4/8) rectangle (\n*8/8, \n*5/8) node[midway] {};

            \draw[fill=\lrcolor, dashed] (\n*0/8, \n*3/8) rectangle (\n*1/8, \n*4/8) node[midway] {};
            \draw[fill=\lrcolor, dashed] (\n*1/8, \n*3/8) rectangle (\n*2/8, \n*4/8) node[midway] {};
            \draw[fill=\lrcolor, dashed] (\n*2/8, \n*3/8) rectangle (\n*3/8, \n*4/8) node[midway] {};
            \draw[fill=\lrcolor, dashed] (\n*6/8, \n*3/8) rectangle (\n*7/8, \n*4/8) node[midway] {};
            \draw[fill=\lrcolor, dashed] (\n*7/8, \n*3/8) rectangle (\n*8/8, \n*4/8) node[midway] {};

            \draw[fill=\lrcolor, dashed] (\n*0/8, \n*2/8) rectangle (\n*1/8, \n*3/8) node[midway] {};
            \draw[fill=\lrcolor, dashed] (\n*1/8, \n*2/8) rectangle (\n*2/8, \n*3/8) node[midway] {};
            \draw[fill=\lrcolor, dashed] (\n*2/8, \n*2/8) rectangle (\n*3/8, \n*3/8) node[midway] {};
            \draw[fill=\lrcolor, dashed] (\n*3/8, \n*2/8) rectangle (\n*4/8, \n*3/8) node[midway] {};
            \draw[fill=\lrcolor, dashed] (\n*7/8, \n*2/8) rectangle (\n*8/8, \n*3/8) node[midway] {};

            \draw[fill=\lrcolor, dashed] (\n*0/8, \n*1/8) rectangle (\n*1/8, \n*2/8) node[midway] {};
            \draw[fill=\lrcolor, dashed] (\n*1/8, \n*1/8) rectangle (\n*2/8, \n*2/8) node[midway] {};
            \draw[fill=\lrcolor, dashed] (\n*2/8, \n*1/8) rectangle (\n*3/8, \n*2/8) node[midway] {};
            \draw[fill=\lrcolor, dashed] (\n*3/8, \n*1/8) rectangle (\n*4/8, \n*2/8) node[midway] {};
            \draw[fill=\lrcolor, dashed] (\n*4/8, \n*1/8) rectangle (\n*5/8, \n*2/8) node[midway] {};

            \draw[fill=\lrcolor, dashed] (\n*0/8, \n*0/8) rectangle (\n*1/8, \n*1/8) node[midway] {};
            \draw[fill=\lrcolor, dashed] (\n*1/8, \n*0/8) rectangle (\n*2/8, \n*1/8) node[midway] {};
            \draw[fill=\lrcolor, dashed] (\n*2/8, \n*0/8) rectangle (\n*3/8, \n*1/8) node[midway] {};
            \draw[fill=\lrcolor, dashed] (\n*3/8, \n*0/8) rectangle (\n*4/8, \n*1/8) node[midway] {};
            \draw[fill=\lrcolor, dashed] (\n*4/8, \n*0/8) rectangle (\n*5/8, \n*1/8) node[midway] {};
            \draw[fill=\lrcolor, dashed] (\n*5/8, \n*0/8) rectangle (\n*6/8, \n*1/8) node[midway] {};

            \draw[fill=\frcolor, dashed] (\n*1/8, \n*7/8) rectangle (\n*2/8, \n*8/8) node[midway] {};
            \draw[fill=\frcolor, dashed] (\n*2/8, \n*6/8) rectangle (\n*3/8, \n*7/8) node[midway] {};
            \draw[fill=\frcolor] (\n*3/8, \n*5/8) rectangle (\n*4/8, \n*6/8) node[midway] {};
            \draw[fill=\frcolor, dashed] (\n*4/8, \n*4/8) rectangle (\n*5/8, \n*5/8) node[midway] {};
            \draw[fill=\frcolor, dashed] (\n*5/8, \n*3/8) rectangle (\n*6/8, \n*4/8) node[midway] {};
            \draw[fill=\frcolor, dashed] (\n*6/8, \n*2/8) rectangle (\n*7/8, \n*3/8) node[midway] {};
            \draw[fill=\frcolor, dashed] (\n*7/8, \n*1/8) rectangle (\n*8/8, \n*2/8) node[midway] {};

            \draw[fill=\frcolor, dashed] (\n*0/8, \n*6/8) rectangle (\n*1/8, \n*7/8) node[midway] {};
            \draw[fill=\frcolor] (\n*1/8, \n*5/8) rectangle (\n*2/8, \n*6/8) node[midway] {};
            \draw[fill=\frcolor, dashed] (\n*2/8, \n*4/8) rectangle (\n*3/8, \n*5/8) node[midway] {};
            \draw[fill=\frcolor, dashed] (\n*3/8, \n*3/8) rectangle (\n*4/8, \n*4/8) node[midway] {};
            \draw[fill=\frcolor, dashed] (\n*4/8, \n*2/8) rectangle (\n*5/8, \n*3/8) node[midway] {};
            \draw[fill=\frcolor, dashed] (\n*5/8, \n*1/8) rectangle (\n*6/8, \n*2/8) node[midway] {};
            \draw[fill=\frcolor, dashed] (\n*6/8, \n*0/8) rectangle (\n*7/8, \n*1/8) node[midway] {};

            \draw[fill=\frcolor, dashed] (\n*0/8, \n*7/8) rectangle (\n*1/8, \n*8/8) node[midway] {};
            \draw[fill=\frcolor, dashed] (\n*1/8, \n*6/8) rectangle (\n*2/8, \n*7/8) node[midway] {};
            \draw[fill=\frcolor] (\n*2/8, \n*5/8) rectangle (\n*3/8, \n*6/8) node[midway] {};
            \draw[fill=\frcolor, dashed] (\n*3/8, \n*4/8) rectangle (\n*4/8, \n*5/8) node[midway] {};
            \draw[fill=\frcolor, dashed] (\n*4/8, \n*3/8) rectangle (\n*5/8, \n*4/8) node[midway] {};
            \draw[fill=\frcolor, dashed] (\n*5/8, \n*2/8) rectangle (\n*6/8, \n*3/8) node[midway] {};
            \draw[fill=\frcolor, dashed] (\n*6/8, \n*1/8) rectangle (\n*7/8, \n*2/8) node[midway] {};
            \draw[fill=\frcolor, dashed] (\n*7/8, \n*0/8) rectangle (\n*8/8, \n*1/8) node[midway] {};

            \draw[fill=\frcolor] (\n*2/8, \n*5/8) rectangle (\n*3/8, \n*6/8) node[midway] {};
            \draw[fill=\frcolor] (\n*1/8, \n*5/8) rectangle (\n*2/8, \n*6/8) node[midway] {};
            \draw[fill=\unifrowscolor] (\n*0/8, \n*5/8) rectangle (\n*1/8, \n*6/8) node[midway] {};
            \draw[fill=\frcolor] (\n*3/8, \n*5/8) rectangle (\n*4/8, \n*6/8) node[midway] {};
            \draw[fill=\unifrowscolor] (\n*4/8, \n*5/8) rectangle (\n*5/8, \n*6/8) node[midway] {};
            \draw[fill=\unifrowscolor] (\n*5/8, \n*5/8) rectangle (\n*6/8, \n*6/8) node[midway] {};
            \draw[fill=\unifrowscolor] (\n*6/8, \n*5/8) rectangle (\n*7/8, \n*6/8) node[midway] {};
            \draw[fill=\unifrowscolor] (\n*7/8, \n*5/8) rectangle (\n*8/8, \n*6/8) node[midway] {};
        \end{scope}
    \end{tikzpicture}
 &  \begin{tikzpicture}[scale=0.8, line width=1pt]
            \def\n{4cm}
            \draw[fill=\unifrowscolor] (\n*0/8, \n*7/8) rectangle (\n*0.66/8, \n*8/8) node[midway] {};
            \draw[fill=white] (\n*0/8, \n*6/8) rectangle (\n*0.66/8, \n*7/8) node[midway] {};
            \draw[fill=white] (\n*0/8, \n*5/8) rectangle (\n*0.66/8, \n*6/8) node[midway] {};
            \draw[fill=white] (\n*0/8, \n*4/8) rectangle (\n*0.66/8, \n*5/8) node[midway] {};
            \draw[fill=\unifrowscolor] (\n*0/8, \n*3/8) rectangle (\n*0.66/8, \n*4/8) node[midway] {};
            \draw[fill=\unifrowscolor] (\n*0/8, \n*2/8) rectangle (\n*0.66/8, \n*3/8) node[midway] {};
            \draw[fill=\unifrowscolor] (\n*0/8, \n*1/8) rectangle (\n*0.66/8, \n*2/8) node[midway] {};
            \draw[fill=\unifrowscolor] (\n*0/8, \n*0/8) rectangle (\n*0.66/8, \n*1/8) node[midway] {};
\end{tikzpicture}
    \end{array},
\end{align*}
\end{minipage}}
\end{center}

\noindent noting that the white blocks of $\momega \mP^{(3)}$ are filled with zeros. 
We also note that the blue blocks of $\momega \mP^{(3)}$ contain standard Gaussian entries
because (1) the distribution of Gaussian matrices is invariant under unitary transformations and (2) the matrix $\mP^{(3)}$ is computed independently of the blue blocks of $\momega$.
Then we can compute $
    \mU_{3} = \col(\mY_{3} \mP^{(3)}, k)
$
with the usual probabilistic guarantees (cf.~Section \ref{sec:preliminaries_randomized} and \cite{halko2011finding, martinsson2020randomized}).

In general, block nullification computes basis matrices $\mU, \mV$ according to Algorithm~\ref{alg:stepI_bn} to accomplish step~\ref{comp1} of uniform BLR compression.
Quickly summarizing, we first draw independent Gaussian matrices $\momega$ and $\mpsi$ to form random sketches $\mY$ and $\mZ$ (lines 1-3).
For any block $i=1,\ldots,b$, we define $\momega^{(i)}$ as the rows of $\momega$ indexed by $\{I_j\}_{j \in \mathcal{N}_i}$, or all $I_j$ such that $j$ is a neighbor of block $i$.
We then compute $r$ orthonormal vectors in the null space of $\momega^{(i)}$, which comprise the columns of  $\mP^{(i)}$.
The matrix $\mY_i = \mY(I_i,:)$ is right-multiplied by $\mP^{(i)}$ to compute $\mU_i$, whose $k$ columns form an approximate basis for the column space of $\ma(I_i,:)$ excluding inadmissible blocks (lines 5-6).
Analogously, for each block $i$, we compute the basis matrix $\mV_{i}$ whose column space approximates the column space of $\ma^*(I_i,:)$ excluding inadmissible blocks (lines 7-8).
We discuss the asymptotic complexity of Algorithm~\ref{alg:stepI_bn} in the next section.

\begin{algorithm}[h]
    \caption{Block Nullification for Basis Construction}\label{alg:stepI_bn} \textbf{Call:} $[\mU, \mV] = \texttt{bn}(\ma, \ma^*, k)$
    \begin{algorithmic}[1]
        \Require Fast matrix-vector multiplication with uniform BLR $\ma \in \mathbb{R}^{N \times N}$ and $\ma^* \in \mathbb{R}^{N \times N}$, target rank $k$, given $b \times b$ tessellation of $d$-dimensional geometry with max block-size $m$
        \Ensure $\mU, \mV \in \mathbb{R}^{N \times bk}$ in uniform BLR representation of $\ma$ as in (\ref{eq:uniformBLRfull})

        \State Set $r = k+p$ and $s \geq r+3^dm$
        \State Draw independent Gaussian matrices $\momega, \mpsi \in \real^{N \times s}$
        \State Sketch $\mY = \ma \momega$ and $\mZ = \ma^* \mpsi$
        \For{blocks $i =1,\ldots,b$}
            \State Compute $\mP^{(i)} = \nullsp \left ( \momega^{(i)}, \ r \right ) $ 
            \State Compute $\mU_i = \col \left (\mY_i \mP^{(i)}, \ k \right )$
            \State  Compute $\mQ^{(i)} = \nullsp \left ( \mpsi^{(i)}, \ r \right ) $ 
            \State Compute $\mV_{i} = \col \left ( \mZ_{i} \mQ^{(i)}, \ k \right )$
        \EndFor
        \State Set $\mU = \textrm{diag} \left (\mU_1,\ldots,\mU_b \right )$ and $\mV = \textrm{diag} \left (\mV_1,\ldots,\mV_b \right ) $
    \end{algorithmic}
\end{algorithm}

\subsection{Asymptotic complexity of block nullification}
\label{sec:blocknull_complexity}
We now analyze the asymptotic complexity of Algorithm~\ref{alg:stepI_bn}. Let $r = k + p$, where $k$ is the block rank and $p$ is the oversampling parameter, and let $s = 3^d m + r$, with $m$ denoting the block size. For the purpose of generality, we do not assume that $\mtx A$ is self-adjoint. However, if that is the case, there are savings of a factor of 2.

\begin{itemize}
\item \textbf{Gaussian matrix generation} (lines 2-3). To generate $\momega$ and $\mpsi$ we draw $2Ns$ values from the standard Gaussian distribution. The cost of this step is $2Ns \times T_{\textup{rand}}$, where $T_{\textup{rand}}$ represents the time to draw a single value.
\item \textbf{Matrix-vector products} (line 3). Forming $\mY = \ma \momega$ and $\mZ = \ma^* \mpsi$ requires $s$ matrix-vector multiplications for $\ma$ and $\ma^*$.  
   This contributes $2s \times T_{\textup{mult}}$ to the overall complexity, where $T_{\textup{mult}}$ is the cost of applying $\ma$ or $\ma^*$ to a vector.
   \textcolor{black}{We note that in the problems we consider (cf. Section~\ref{sec:num_exp}), we have fast matrix-vector products so that $T_{mult}$ is $N \log N$.}

\item \textbf{Null-space bases} (lines 5, 7).  
   For block $i = 1, \ldots, b$, we compute $\mP^{(i)} = \nullsp(\momega^{(i)}, r)$ and $\mQ^{(i)} = \nullsp(\mpsi^{(i)}, r)$ with matrices of size at most $3^dm \times s$. Using Householder QR \cite[Table~C.2]{Higham2008} for $\nullsp$, the cost of each is roughly
   \(
   \mathcal{O}\left(3^{3d} m^3 \right).
   \)
   For $b$ blocks, this step adds
   $
   \mathcal{O}\left(b \times 3^{3d} m^3 \right) \times T_{\textup{flop}},$
   where $T_{\textup{flop}}$ is the cost of one floating-point arithmetic operation.
   \item \textbf{Column basis extraction} (lines 6, 8). 
   Computing $\mU_i = \col(\mY_i \mP^{(i)}, k)$ and $\mV_i = \col(\mZ_i \mQ^{(i)}, k)$ involves matrix multiplications of dimensions $3^dm \times r$ by $r \times k$ for each block.  
   The cost for all blocks is
   $
   \mathcal{O}\left(b \times [3^dmr^2]\right) \times T_{\textup{flop}}.$
   \item  \textbf{Reconstruction.}  
   Once $\mU$ and $\mV$ are computed, determining the matrix $\mB$ for the full matrix reconstruction as in (\ref{eq:uniformBLRfull}) involves additional matrix-vector products. The cost of this step is exactly
   $
   (3^dm + kb) \times T_{\textup{mult}}.
   $
\end{itemize}

Combining the contributions from all steps and using $N = mb$ gives total cost
\begin{align*}
2N(3^dm + r) \times T_{\textup{rand}} &+ 
\mathcal{O}\left(N \times 3^{3d}m^2\right) \times T_{\textup{flop}} \\
&+ (3^{(d+1)}m + 2r + kN/m) \times T_{\textup{mult}}.
\end{align*}
\textcolor{black}{The key takeaway here is that both the number of samples needed to compute the basis matrices with block nullification and the null space computations depend on the block size $m$, and consequently, the problem size $N$, which is computationally undesirable.
Our proposed method of tagging avoids these dependencies on the problem size, as discussed in the next section.}


\section{Tagging in Uniform BLR Matrix Compression}
\label{sec:tagging}

We now introduce a black-box randomized method to compress strongly admissible uniform BLR matrices which we call tagging, the main contribution of our manuscript. 
As in Section~\ref{sec:blocknull}, we begin with an illustrative example before generalizing the tagging method for $d$-dimensional problem geometries and summarizing its asymptotic complexity. 

\subsection{Tagging for strongly-admissible uniform BLR matrices}
\label{sec:tagging_overview}

Let $\ma$ be an $N \times N$ strongly admissible uniform BLR matrix as in Fig.~\ref{fig:blr2}, tessellated into $b \times b$ blocks each of size $m \times m$ with block rank $k$, allowing for a small amount of oversampling given by $p$, and let $r = k+p$.
Our aim once again is to construct random test matrices $\momega, \mpsi \in \mathbb{R}^{N \times s}$ with $s = \mathcal{O}(k)$ so that the sketches $\mY = \ma \momega$ and $\mZ = \ma^* \mpsi$ taken a priori can be used to compress the admissible blocks in every block-row and block-column. 

To illustrate, suppose as in Section~\ref{sec:blocknull_overview} that we want to compute the basis matrix $\mU_{3} \in \mathbb{R}^{m \times k}$ from Fig.~\ref{fig:unifbasis}.
We first introduce the $8 \times 4$ \textit{tagging matrix} 
\begin{align}
\label{eq:tagmat}
            \mT
    =
    \begin{bmatrix}
        t_{1, 1} & t_{1, 2} & t_{1, 3} & t_{1, 4} \\
        t_{2, 1} & t_{2, 2} & t_{2, 3} & t_{2, 4} \\
        \vdots & \vdots & \vdots & \vdots \\
        t_{8, 1} & t_{8, 2} & t_{8, 3} & t_{8, 4} \\
    \end{bmatrix},
\end{align}
where the entries $t_{i,j}$ will be made explicit in Section~\ref{sec:tagmat_selection}; for now, we assume standard Gaussian entries.
The number of rows of $\mT$ equals the number of matrix blocks $b$.
The number of columns of $\mT$ is one more than the maximal number of neighbors within a given block-row/column, which is $4$ for the strongly admissible 1-D matrix in Fig.~\ref{fig:blr2}.
In general, for a $d$-dimensional problem geometry, the number of columns of $\mT$ is $3^d+1$.

We next define the \textit{extended random test matrix} $\momega \in \mathbb{R}^{N \times 4r}$ in terms of random test matrices $\momega_j \in \mathbb{R}^{N \times r}$ for $j=1,\ldots,4$, given by
\begin{align}
\label{eq:extend_test_mat}
        \underset{N \times 4r}{\momega} := \begin{bmatrix} {\momega_1} & \momega_2 & \momega_3 & \momega_4  \end{bmatrix} = \begin{bmatrix}
            t_{1,1} \mG_{1} & t_{1,2} \mG_{1} & t_{1,3} \mG_{1} & t_{1,4} \mG_{1} \\
            t_{2,1} \mG_{2} & t_{2,2} \mG_{2} & t_{2,3} \mG_{2} & t_{2,4} \mG_{2} \\
            t_{3,1} \mG_{3} & t_{3,2} \mG_{3} & t_{3,3} \mG_{3} & t_{3,4} \mG_{3} \\
            t_{4,1} \mG_{4} & t_{4,2} \mG_{4} & t_{4,3} \mG_{4} & t_{4,4} \mG_{4} \\
            \vdots & \vdots & \vdots & \vdots \\
            t_{8,1} \mG_{8} & t_{8,2} \mG_{8} & t_{8,3} \mG_{8} & t_{8,4} \mG_{8} 
        \end{bmatrix},
\end{align}
where each $\mG_i \in \mathbb{R}^{m \times r}$, $i = 1,\ldots,8$, is a Gaussian random matrix, weighted by entry $t_{i,j}$ of the tagging matrix $\mT$ to form $\momega_j$, $j = 1,\ldots,4$.
Note that we again assume $N = bm$ for notational convenience.
We form the sketch matrix $\mY = \ma \momega \in \mathbb{R}^{N \times 4r}$, partitioned into $N \times r$ block-columns commensurately with (\ref{eq:extend_test_mat}) so that
\begin{align*} \underset{N \times 4r}{\mY} = \begin{bmatrix}
        \mY_1 & \mY_2 & \mY_3 & \mY_4
    \end{bmatrix} = \begin{bmatrix}
        \ma \momega_1 & \ma \momega_2 & \ma \momega_3 & \ma \momega_4
    \end{bmatrix} = \ma \momega. \end{align*}
To compute $\mU_{3} \in \mathbb{R}^{m \times k}$ as in Fig.~\ref{fig:unifbasis}, we exclude contributions from inadmissible blocks by computing a (nonzero) vector $\vz^{(3)} = \begin{bmatrix}
            z_1^{(3)} & z_2^{(3)} & z_3^{(3)} & z_4^{(3)}
        \end{bmatrix}^T$ 
so that
\begin{align}
    \label{eq:tagsubmat}
    \vz^{(3)} = \nullsp \left ( \mT^{(3)} \right ) :=
    \nullsp \left ( \begin{bmatrix}
        t_{2, 1} & t_{2, 2} & t_{2, 3} & t_{2, 4} \\
        t_{3, 1} & t_{3, 2} & t_{3, 3} & t_{3, 4} \\
        t_{4, 1} & t_{4, 2} & t_{4, 3} & t_{4, 4}
    \end{bmatrix} \right ).
\end{align}
Note that this submatrix $\mT^{(3)}$ of $\mT$ comprises the rows that correspond to the neighbor list $\mathcal{N}_3 = \{2,3,4\}$ of box 3 (Fig.~\ref{fig:blr2}).
Now consider the weighted sum
    \begin{align*}
        z_1^{(3)} \momega_1 + z_2^{(3)} \momega_2 + z_3^{(3)} \momega_3 + z_4^{(3)} \momega_4 = \begin{bmatrix}
           (z_1^{(3)} t_{1,1} + z_2^{(3)} t_{1,2} + z_3^{(3)} t_{1,3} + z_4^{(3)} t_{1,4}) \mG_{1}\\
           (z_1^{(3)} t_{2,1} + z_2^{(3)} t_{2,2} + z_3^{(3)} t_{2,3} + z_4^{(3)} t_{2,4})  \mG_{2}\\
            (z_1^{(3)} t_{3,1} + z_2^{(3)} t_{3,2} + z_3^{(3)} t_{3,3} + z_4^{(3)} t_{3,4})  \mG_{3}  \\
            (z_1^{(3)} t_{4,1} + z_2^{(3)} t_{4,2} + z_3^{(3)} t_{4,3} + z_4^{(3)} t_{4,4}) \mG_{4}  \\
            (z_1^{(3)} t_{5,1} + z_2^{(3)} t_{5,2} + z_3^{(3)} t_{5,3} + z_4^{(3)} t_{5,4}) \mG_{5}  \\
            \vdots \\
            (z_1^{(3)} t_{8,1} + z_2^{(3)} t_{8,2} + z_3^{(3)} t_{8,3} + z_4^{(3)} t_{8,4})  \mG_{8} 
        \end{bmatrix}.
    \end{align*}
By construction, it simplifies to
\begin{align}
\label{eq:projected_tags}
z_1^{(3)} \momega_1 + z_2^{(3)} \momega_2 + z_3^{(3)} \momega_3 + z_4^{(3)} \momega_4 = 
\begin{bmatrix}
           (z_1^{(3)} t_{1,1} + z_2^{(3)} t_{1,2} + z_3^{(3)} t_{1,3} + z_4^{(3)} t_{1,4}) \mG_{1}\\
           \mathbf{0} \\
            \mathbf{0}  \\
            \mathbf{0} \\
            (z_1^{(3)} t_{5,1} + z_2^{(3)} t_{5,2} + z_3^{(3)} t_{5,3} + z_4^{(3)} t_{5,4}) \mG_{5}  \\
            \vdots \\
            (z_1^{(3)} t_{8,1} + z_2^{(3)} t_{8,2} + z_3^{(3)} t_{8,3} + z_4^{(3)} t_{8,4})  \mG_{8} 
        \end{bmatrix},
\end{align}
so that the rows of $\ma(z_1^{(3)} \momega_1 + z_2^{(3)} \momega_2 + z_3^{(3)} \momega_3 + z_4^{(3)} \momega_4) = z_1^{(3)} \mY_1 + z_2^{(3)} \mY_2 + z_3^{(3)} \mY_3 + z_4^{(3)} \mY_4$ corresponding to $I_3$ will contain the desired sample of the admissible blocks, with the contributions from inadmissible blocks now zeroed out.
We can then compute $\mU_3 = \col \left ( (z_1^{(3)} \mY_1 + z_2^{(3)} \mY_2 + z_3^{(3)} \mY_3 + z_4^{(3)} \mY_4)(I_3,:),k \right )$ as in Section~\ref{sec:preliminaries_randomized}.

In general, the tagging method computes basis matrices $\mU$ and $\mV$ according to Algorithm~\ref{alg:stepI_tagging}, which we quickly summarize.
We begin by drawing entries of the tagging matrix $\mT \in \mathbb{R}^{b \times (3^d+1)}$, e.g., from a standard Gaussian distribution (line 2).
We next form the random sketches 
\begin{align}
    \underset{N \times (3^d+1)r}{\mY} &= \begin{bmatrix}
         \mY_1 & \mY_2 & \ldots & \mY_{3^d+1}
    \end{bmatrix} = \begin{bmatrix}
        \ma \momega_1 & \ma \momega_2 & \ldots & \ma \momega_{3^d+1}
    \end{bmatrix} = \ma \momega, \\
    \underset{N \times (3^d+1)r}{\mZ} &= \begin{bmatrix}
        \mZ_1 & \mZ_2 & \ldots & \mZ_{3^d+1}
    \end{bmatrix} = \begin{bmatrix}
        \ma^* \mpsi_1 & \ma^* \mpsi_2 & \ldots & \ma^* \mpsi_{3^d+1}
    \end{bmatrix} = \ma^* \mpsi.
\end{align}
where each $\momega_j$, $j=1,\ldots,3^d+1$, comprises $b$ block-rows of independent Gaussian matrices $\mG_i \in \mathbb{R}^{m \times r}$, $i=1,\ldots,b$, weighted by tagging entry $t_{i,j}$ as in (\ref{eq:extend_test_mat}), and similarly for each $\mpsi_j$ with $\mtx{H}_i$ (lines 3-8).
To compute $\mU_{i}$ or $\mV_{i}$ for block $i$, we focus on the submatrix $\mT^{(i)}$ comprising rows of $\mT$ indexed by $\mathcal{N}_{i}$, the neighbor list of box $i$, so that $|\mathcal{N}_{i}| \leq 3^d$.
We then compute an orthonormal tagging vector $\vz^{(i)} = \begin{bmatrix} z^{(i)}_1 z^{(i)}_2 \ldots z^{(i)}_{3^d+1} \end{bmatrix}^T$ in the nontrivial null space of $\mT^{(i)}$ (line 10).
Finally, we compute $\mU_{i}, \mV_{i} \in \mathbb{R}^{m \times k}$ (lines 11-12) whose columns are an ON-basis approximating the column spaces of $\ma(I_i,:)$ and $\ma^*(I_i,:)$, excluding contributions from the inadmissible blocks.

\begin{algorithm}[h]
    \caption{Tagging for Basis Construction}\label{alg:stepI_tagging} \textbf{Call:} $[\mU,\mV] = \texttt{tag}(\ma, \ma^*, k)$
    \begin{algorithmic}[1]
        \Require Fast matrix-vector multiplication with uniform BLR $\ma \in \mathbb{R}^{N \times N}$ and $\ma^* \in \mathbb{R}^{N \times N}$, target rank $k$, given $b \times b$ tessellation of $d$-dimensional geometry with max block-size $m$ 
        \Ensure $\mU, \mV \in \mathbb{R}^{N \times bk}$ in uniform BLR representation of $\ma$ as in (\ref{eq:uniformBLRfull})

        \State Set $r = k+p$ and $\momega, \mpsi = [ \ ]$.
        \State Form tagging matrix $\mT \in \mathbb{R}^{b \times (3^d+1)}$  \LineComment{e.g. Gaussian $\mT$; see Section~\ref{sec:tagmat_selection}}
        \For{blocks $i=1,\ldots,b$}
            \State Draw independent Gaussian matrices $\mtx{G}_i, \mtx{H}_i \in \real^{m \times r}$
            \State Update $\momega = \begin{bmatrix} \momega & ; & \begin{bmatrix} t_{i,1} \mtx G_i & \ldots & t_{i,3^d+1} \mtx G_i \end{bmatrix} \end{bmatrix}$ \LineComment{Append $i^{th}$ row to $\momega$}
            \State Update $\mpsi = \begin{bmatrix} \mpsi & ; & \begin{bmatrix} t_{i,1} \mtx H_i & \ldots & t_{i,3^d+1} \mtx H_i \end{bmatrix} \end{bmatrix}$ \LineComment{Append $i^{th}$ row to $\mpsi$}
        \EndFor
        \State Sketch $\mY = \ma \momega$ and $\mZ = \ma^* \mpsi$
        \For{blocks $i=1,\ldots,b$}
            \State Compute $\mtx{z}^{ (i)} = \nullsp \left ( \mT^{(i)} \right )$ 
            \State Compute $\mU_i = \col \left ( \sum_{j=1}^{3^d+1} z_j^{(i)} \mY_j(I_i, :), \ k \right )$
            \State Compute $\mV_i = \col \left (\sum_{j=1}^{3^d+1} z_j^{(i)} \mZ_j(I_i, :), \ k \right )$
        \EndFor
        \State $\mU = \textrm{diag} \left (\mU_1,\ldots,\mU_b \right )$ and $\mV = \textrm{diag} \left (\mV_1,\ldots,\mV_b \right ) $
    \end{algorithmic}
\end{algorithm}

\subsection{Asymptotic Complexity of Tagging}
\label{sec:tagging_complexity}

We derive the asymptotic complexity of Algorithm~\ref{alg:stepI_tagging} in terms of the problem size $N$, block size $m$, block-rank $k$, oversampling parameter $p$, and problem geometry dimension $d$. Let $r = k + p$ and assume $N = mb$, where $b$ is the number of matrix blocks in each block-row or block-column.
The complexity of the algorithm is broken into the following components:

\begin{itemize}
    \item \textbf{Tagging matrix formation} (line 2).
    The tagging matrix $\mT \in \mathbb{R}^{b \times (3^d+1)}$ requires $(3^d+1)N/m$ elements to be sampled from a Gaussian distribution. The total complexity is
    $
    (3^d+1)N/m \times T_{\textup{rand}}.
    $

    \item \textbf{Random test matrix formation} (lines 3--7).
    Forming $\momega, \mpsi \in \mathbb{R}^{N \times (3^d+1)r}$ first requires sampling $\mtx G_i, \mtx H_i$ for each subblock and weighting each by the relevant tag.
    This contributes $2Nr \times T_{\textup{rand}}$, as each Gaussian matrix $\mtx{G}_i$ and $\mtx{H}_i$ is reused.
    Forming $\momega,\mpsi$ then requires $2(3^d+1)Nr$ floating-point operations. The total complexity is $
    2 N r \times T_{\textup{rand}} + \mathcal{O}(3^d N r) \times T_{\textup{flop}}.$
    \item \textbf{Matrix-vector products} (line 8). 
    Computing $\mY = \ma \momega$ and $\mZ = \ma^* \mpsi$ involves $2(3^d+1)r$ matrix-vector products, with complexity $
    2(3^d+1)r \times T_{\textup{mult}}.
    $
    \textcolor{black}{We note that in the problems we consider (cf. Section~\ref{sec:num_exp}), we have fast matrix-vector products so that $T_{mult}$ is $N \log N$.}

    \item \textbf{Null space bases} (line 10).
    For each block, the null space computation involves matrices of size at most $3^d \times (3^d+1)$. Using Householder QR, the complexity is $
    \mathcal{O}\left(3^{3d} b \right) \times T_{\textup{flop}}.$

    \item \textbf{Column basis extraction} (lines 11, 12). 
    Forming $\mU$ and $\mV$ involves multiplying the relevant subblocks by the null space vector, then computing an ON-basis.
    This involves $(3^d+1)mr$ floating-point operations per block to scale submatrices of $\mY$ and $\mZ$ and $\mathcal{O}(rmk)$ floating-point operations for $\mU_i$ and $\mV_i$.
    The total complexity is $
    \mathcal{O}(3^dNr + Nrk) \times T_{\textup{flop}}.$

    \item \textbf{Reconstruction.} 
    Reconstructing the uniform BLR matrix in (\ref{eq:uniformBLRfull}) using $\mU$ and $\mV$ contributes an additional $
    (3^dm + kb) \times T_{\textup{mult}}$.
\end{itemize}
Combining the above contributions, the total complexity is
\begin{align*}
\left((3^d+1)N/m + 2Nr\right) \times T_{\textup{rand}}
 &+ \mathcal{O}\left(3^dNr + Nrk + 3^{3d}N/{m}\right) \times T_{\textup{flop}} \\
 &+ (2 (3^d+1) r + 3^d m + kb) \times T_{\textup{mult}}.
\end{align*}
Compared to block nullification in Section \ref{sec:blocknull_complexity},  far fewer samples are needed to construct $\mtx U$ and $\mtx V$. Block nullification requires $2\times 3^{d} m$ samples for basis construction in Section \ref{sec:blocknull_complexity}, whereas tagging only
requires $2 \times 3^{d+1} r$ samples. The cost of reconstructing the uniform BLR matrix, however,
dominates the asymptotic complexity of $T_{\textup{mult}}$ for both methods.
The key advantage of tagging is the substantially reduced (linear) cost of post-processing the test and sketch matrices, as opposed to block-nullification, where post-processing scales linearly with $N$ and quadratically with $m$.

\section{Selecting the Tagging Matrix}
\label{sec:tagmat_selection}

In Section~\ref{sec:tagging_overview}, we treated the tagging matrix as having standard Gaussian entries.
However, standard Gaussian tagging matrix entries do not guarantee Gaussian samples of the input matrix.
As such, we seek to address the following questions. \textit{Does there exist an optimal tagging matrix? If not, can we still obtain a high-quality randomized embedding?}

To answer these questions, we first discuss our criteria for tagging matrix optimality in Section~\ref{sec:tagmat_optimality}.
We then present a conjecture on the existence of optimal tagging matrices in Section~\ref{sec:tagmat_existence} which takes an algebraic-geometric perspective.
We finish the section with a highly efficient alternative strategy to determine tagging matrices that perform well empirically despite their sub-optimality.

\subsection{On the optimality of tagging matrices: Projected tags and aspect ratios}
\label{sec:tagmat_optimality}

The main issue that we need to address in tagging matrix selection concerns the \textit{projected tags}, the nonzero non-uniform weights on each Gaussian matrix in (\ref{eq:projected_tags}).
These are given by $t_{j,1} z_1^{(i)} + t_{j,2} z_2^{(i)} + \ldots t_{j,\ell} z_{\ell}^{(i)}$ with $1 \leq j \leq b$ for each box $i=1,\ldots,b$. 
Here $\ell = 3^d+1$ and $\vz^{(i)} = \begin{bmatrix} z^{(i)}_1 & z^{(i)}_2 & \ldots & z^{(i)}_{\ell} \end{bmatrix}^T = \nullsp (\mT^{(i)})$ for $\mT^{(i)} = \mT(\mathcal{N}_i,:)$ where $\mathcal{N}_i$ is the neighbor list of box $i$.
The projected tags for $\mG_j$ are 0 by construction for any $j \in \mathcal{N}_i$.
However, the nonzero projected tags, corresponding to the far-field $\mathcal{F}_i = [b] \setminus \mathcal{N}_i$ of box $i$, non-uniformly weight each Gaussian matrix, resulting in non-uniformly weighted randomized samples of blocks within the same block-row or block-column of $\ma$.
Each of the projected tags should ideally be equal in magnitude.
\textcolor{black}{Otherwise, for an adversarially-constructed matrix, blocks that do not contribute much to the basis of the row or column space could be more heavily-weighted than blocks that do; we outline methods in the remainder of this section that address this issue.}

To this end, we define the \textit{aspect ratio} $\rho^{(i)}$ for block-row/column $i=1,\ldots,b$ as the largest-magnitude to the smallest-magnitude nonzero projected tag:
\begin{align}
    \label{eq:aspect_ratio}
    \rho^{(i)} = \frac{\max\limits_{j \in \mathcal{F}_i} \vert t_{j,1}z_1^{(i)} + t_{j,2}z_2^{(i)} + \ldots + t_{j,\ell} z^{(i)}_{\ell} \vert }{\min\limits_{j \in \mathcal{F}_i} |t_{j,1}z_1^{(i)} + t_{j,2}z_2^{(i)} + \ldots + t_{j,\ell} z^{(i)}_{\ell} |} \geq 1.
\end{align}
An optimal tagging matrix then minimizes $\rho^{(i)}$ for each $i=1,\ldots,b$.

In the following subsections, we examine the tagging matrix optimality problem through two different lenses.
The first relies on ideas from algebraic geometry to  determine an optimal tagging matrix. 
The second offers a numerical short-cut via null space vectors that minimize aspect ratios through a fast optimization scheme.

\subsection{On the existence of optimal tagging matrices: Pl\"ucker coordinates}
\label{sec:tagmat_existence}

To express our conjecture on optimal tagging matrices, we draw a connection between tagging matrices and projective varieties through \textit{Pl{\"u}cker coordinates}.
We return to our example of a uniform BLR matrix from Fig.~\ref{fig:blr2} for an intuitive introduction to the Pl{\"u}cker embedding that gives rise to Pl{\"u}cker coordinates.
We then hypothesize that optimal tagging matrices may be found through a hybrid numeric-symbolic approach based on Pl{\"u}cker coordinates, which is currently out of reach.

\subsubsection{An illustrative example}
Consider $\ma \in \mathbb{R}^{N \times N}$ from Fig.~\ref{fig:blr2}. 
Note that for blocks $i = 2,\ldots,7$\footnote{We treat the ``extremal'' blocks $i=1$ and $i=8$ (in general, blocks with fewer than the maximal number of neighbors $3^d$) at the end of the section.}, we can apply Cramer's Rule to find $\vz^{(i)} = \nullsp(\mT^{(i)})$, e.g. the $i^{th}$ coordinate $z_i^{(3)}$ of $\vz^{(3)}$ is the determinant of $\mtx{T}^{(3)}$ without the $i^{th}$ column:
\begin{align}
\begin{split}
\label{eq:cramers_rule}
    z_1^{(3)} &= \det \left ( \begin{bmatrix}
        t_{2,2} & t_{2,3} & t_{2,4} \\
        t_{3,2} & t_{3,3} & t_{3,4} \\
        t_{4,2} & t_{4,3} & t_{4,4}
    \end{bmatrix} \right ), \qquad z_2^{(3)} = \det \left ( \begin{bmatrix}
        t_{2,1} & t_{2,3} & t_{2,4} \\
        t_{3,1} & t_{3,3} & t_{3,4} \\
        t_{4,1} & t_{4,3} & t_{4,4}
    \end{bmatrix} \right ), \\
    z_3^{(3)} &= \det \left ( \begin{bmatrix}
        t_{2,1} & t_{2,2} & t_{2,4} \\
        t_{3,1} & t_{3,2} & t_{3,4} \\
        t_{4,1} & t_{4,2} & t_{4,4}
    \end{bmatrix} \right ), \qquad z_4^{(3)} = \det \left ( \begin{bmatrix}
        t_{2,1} & t_{2,2} & t_{2,3} \\
        t_{3,1} & t_{3,2} & t_{3,3} \\
        t_{4,1} & t_{4,2} & t_{4,3}
    \end{bmatrix} \right ).
\end{split}
\end{align}
Then the projected tag $\mT \vz^{(i)}$ contains all $4 \times 4$ determinants of $\mT$ formed from the three rows of $\mT^{(i)}$ plus one remaining row of $\mT$, e.g. for $\vz^{(3)}$ using (\ref{eq:cramers_rule}),
\begin{align}
\label{eq:4by4}
\begin{split}
    \mT \vz^{(3)} &= \begin{bmatrix} z_1^{(3)} t_{1,1} + z_2^{(3)} t_{1,2} + z_3^{(3)} t_{1,3} + z_4^{(3)} t_{1,4} \\
    \vdots \\ 
    z_1^{(3)} t_{8,1} + z_2^{(3)} t_{8,2} + z_3^{(3)} t_{8,3} + z_4^{(3)} t_{8,4} 
    \end{bmatrix} = \begin{bmatrix}
        \det \left ( \begin{bmatrix}
            t_{1,1} & t_{1,2} & t_{1,3} & t_{1,4} \\
            t_{2,1} & t_{2,2} & t_{2,3} & t_{2,4} \\
            t_{3,1} & t_{3,2} & t_{3,3} & t_{3,4} \\
            t_{4,1} & t_{4,2} & t_{4,3} & t_{4,4} 
        \end{bmatrix} \right ) \\
        \vdots \\
        \det \left ( \begin{bmatrix}
            t_{2,1} & t_{2,2} & t_{2,3} & t_{2,4} \\
            t_{3,1} & t_{3,2} & t_{3,3} & t_{3,4} \\
            t_{4,1} & t_{4,2} & t_{4,3} & t_{4,4} \\
            t_{8,1} & t_{8,2} & t_{8,3} & t_{8,4} \\
        \end{bmatrix} \right ) 
    \end{bmatrix} .
\end{split}
\end{align}
When $\mT$ has full rank, these $4 \times 4$ determinants form a subset of the \textit{Pl{\"u}cker relations}, the set of all possible $4 \times 4$ determinants with the rows of $\mT$.
Let $\mathcal{L}$ be a 4-dimensional subspace of $\mathbb{R}^8$ with $\textrm{Col}(\mT) = \mathcal{L}$.
The \textit{Pl{\"u}cker embedding} maps $\mathcal{L}$ to the point in real projective space whose coordinates are all $4 \times 4$ determinants of $\mT \in \mathbb{R}^{8 \times 4}$.
In other words, the Pl{\"u}cker embedding maps the \textit{Grassmannian manifold} $\textrm{Gr}(4,8)$, comprising all $4$-dimensional subspaces of $\mathbb{R}^8$, to \textit{Pl{\"u}cker coordinates} in $\mathbb{P}^{\binom{8}{4}-1}$, as stated below:

\begin{definition}
    The \textbf{Pl{\"u}cker embedding} is the map $\Gr(k,n) \rightarrow \mathbb{P}^{\binom{n}{k} - 1}$ that identifies $\mathcal{V} \in \Gr(k,n)$ with a unique point in  real projective space with coordinates given by all $k \times k$ determinants of a matrix $\mB \in \mathbb{R}^{n \times k}$ satisfying $\mathcal{V} = \textup{Col}(\mB)$, called \textbf{Pl{\"u}cker coordinates}.
\end{definition}

\begin{definition}
    Let $\mathcal{V} \in \Gr(k,n)$ and suppose $\mB \in \mathbb{R}^{n \times k}$ satisfies $\mathcal{V} = \textup{Col}(\mB)$. 
    For any ordered sequence of $k$ row indices $1 \leq i_1 < \ldots < i_k \leq n$ of $\mB$, let $B_{i_1,\ldots,i_k}$ be the determinant of $k\times k$ submatrix $\mB([i_1,\ldots,i_k],:)$, so that set of all Pl{\"u}cker coordinates may be denoted $\{B_{i_1,\ldots,i_k}\}$.
    For any two ordered sequences of row indices
    \begin{align*}
       1 \leq i_1 < i_2 < \ldots < i_{k-1} \leq n, \qquad 1 \leq j_1 < j_2 < \ldots < j_{k+1} \leq n,
    \end{align*}
    the \textbf{Pl{\"u}cker relations} are the following  homogeneous quadratic equations that must hold for all  Pl{\"u}cker coordinates $\{B_{i_1,\ldots,i_k}\}$:
    \begin{align}
        \label{eq:plucker_relations}
        \sum_{\ell = 1}^{k+1} (-1)^{\ell} B_{i_1,\ldots,i_{k-1},j_{\ell}} B_{j_1,\ldots,\widehat{j}_{\ell},\ldots,j_{k+1}} = 0,
     \end{align}
     where $j_1,\ldots,\widehat{j}_{\ell},\ldots,j_{k+1}$ is the sequence $j_1,\ldots,j_{k+1}$ with the term $j_{\ell}$ omitted.
\end{definition}
We can now pose the optimality of tagging matrices in algebraic-geometric terms.

\subsubsection{On the existence of optimal tagging matrices}
\label{sec:plucker_conjecture}

Let $\mT \in \mathbb{C}^{b \times 3^d+1}$ be a matrix of indeterminates $t_{i,j}$.
For any ordered sequence of $3^d+1$ row indices $1 \leq i_1 < \ldots < i_{3^d+1} \leq b$ of $\mT$, let $T_{i_1,\ldots,i_{3^d+1}}$ be the determinant of the $(3^d+1) \times (3^d+1)$ submatrix  $\mT([i_1,\ldots,i_{3^d+1}],:)$; e.g. $T_{2,3,4,8}$ is the $4 \times 4$ determinant in (\ref{eq:4by4}) involving $\mT^{(3)}$ and the last row of $\mT$.
Each $T_{i_1,\ldots,i_{3^d+1}}$ is a degree-$(3^d+1)$ polynomial in the indeterminates $t_{i,j}$, and from the previous section, each is a Pl\"ucker coordinate that must satisfy the Pl\"ucker relations, which are quadratic polynomials in the indeterminates $T_{i_1,\ldots,i_{3^d+1}}$ for all possible ordered sequences $1 \leq i_1 < \ldots < i_{3^d+1} \leq b$.

To determine tagging matrix entries $t_{i,j}$ that minimize $\rho^{(i)}$ for each block $i$, we propose the following approach. 
It is well-known that the set of Pl\"ucker relations is not algebraically independent, cf. \cite[Chapter 14.2]{miller2004combinatorial} and \cite[Appendix C.7]{Harnad_Balogh_2021}.
Thus, the first step is to determine an algebraically independent generating set of  Pl\"ucker relations for $\mT$.
One method is the computation of a Gr\"obner basis\footnote{The computation of a Gr\"obner basis for the quadratic polynomials under consideration is highly nontrivial for problems of this size due to, e.g., intermediate swell in Buchberger's algorithm; see \cite{Demin2023GroebnerjlAP}.} for the ideal of the polynomial ring $\mathbb{C}[\{t_{i,j} \}_{i \in [b], j \in [3^d+1]}]$ generated by all Pl\"ucker relations, cf. \cite{MouChenqi2024DASP}.

{\remark{Another avenue of investigation that bears future consideration involves the so-called ``clusters'' formed by independent Pl\"ucker coordinates \cite{Scott2006}. 
One such cluster is comprised of \textit{rectangular Pl\"ucker coordinates} \cite{Karp2019}, which correspond to the rectangular partitions of a $k \times (n-k)$ unit rectangle and form a generating set for the coordinate ring. 
Rectangular Pl\"ucker coordinates relate to quantum Schubert calculus on the flag variety \cite{Batyrev2000, Eguchi1997, Givental1997} and have an associated Laurent polynomial with certain properties \cite{MARSH2020107027} that may offer another path to an optimal tagging matrix.}}
\vspace{2mm}

Let $\mathcal{R}$ denote a set of algebraically independent Pl\"ucker relations that generate the projective variety defined by them all.
For the purposes of tagging, we are only interested in the nonzero projected tags, which correspond to the far-field of each box.
Let $n$ be the total number of nonzero projected tags, with $n_i$ nonzero projected tags for block $i$, so that $\sum_{i=1}^b n_i = n$.
Denote the $n$ nonzero projected tags by $\{T_1,\ldots,T_n\}$:
\begin{align}
   \{T_1,\ldots,T_n\}:= \left  \{T_{i^{(1)}_1,\ldots,i^{(1)}_{3^d+1}}, \ T_{i^{(2)}_1,\ldots,i^{(2)}_{3^d+1}}, \ldots, \ T_{i^{(n)}_1,\ldots,i^{(n)}_{3^d+1}}  \right \}
\end{align}
for distinct ordered sequences $1 \leq i^{(\ell)}_1 < \ldots < i^{(\ell)}_{3^d+1} \leq b$ for $\ell = 1, \ldots, n$. 

Ideally, for each block $i$, every nonzero projected tag, or nonzero coordinate of $\mT \mtx{z}^{(i)}$, should be equal (cf.~(\ref{eq:aspect_ratio})), which we can enforce numerically via 
\begin{align}
\begin{split}
\label{eq:plucker_opt}
    &\min_{\mtx{T} \in \mathbb{R}^{b \times 3^d+1}} \sum_{i=1}^b  \sum_{\ell_i=1}^{n_i} (T_{\ell_i} - \overline{T}_i)^2 \\
    &\textup{subject to } \ \mathcal{R}
\end{split}
\end{align}
where $\overline{T}_i$ denotes the mean of the $n_i$ nonzero projected tags for block $i$. 

We now note that this approach only holds for blocks with neighbor lists of maximal size $3^d$. 
One workaround is to treat the necessary number of admissible blocks as inadmissible in the block-rows/columns of $\ma$ with fewer than $3^d$ inadmissible blocks, though this increases the overall cost.
More detrimental, though, is the cost of the symbolic computations.
The computations for $\mathcal{R}$ are highly nontrivial even for very small values of $b$, constraints which must hold if the arg min $\mtx{T}^*$ of (\ref{eq:plucker_opt}) satisfies $\textup{Col}(\mtx{T}^*) = \mathcal{V}$ for some $(3^d+1)$-dimensional subspace $\mathcal{V}$ of $\mathbb{R}^{b}$.
Moreover, if the problem size or underlying geometry were to change, these computations would need to be done anew for different values of $b$ or $d$.
We next describe a strictly numerical method to minimize the aspect ratios that performs very well in practice.

\subsection{A practical method for numerical optimization}
\label{sec:tagmat_practical}

Because of the practical difficulties in using the approach of Section~\ref{sec:tagmat_existence}, we  present an alternative to minimize the aspect ratios numerically. 
This approach requires higher-dimensional null spaces of tagging submatrices; as such, we now consider tagging matrices $\mT \in \mathbb{R}^{b \times \ell}$ where $\ell > 3^d+1$ so that every null space has dimension strictly greater than 1.

To minimize the aspect ratios of projected tags efficiently, we consider the following optimization problem over the $(\ell-3^d)$-dimensional null space of $\mT^{(i)} \in \mathbb{R}^{3^d \times \ell}$, rather than over all possible matrix representations $\mT$ of a Grassmannian subspace $\mathcal{V}$ such that $\mathcal{V}=\textup{Col}(\mT)$. 
We now seek a unit vector $\vz^{(i)}$ in the null space of $\mT^{(i)}$ for each block $i = 1,\ldots,b$, which minimizes the ratio of projected tags:
 \begin{align}
    \vz^{(i)} = \arg \min_{\xx \in \textup{Null}(\mT^{(i)})} \frac{\max\limits_{j \in \mathcal{F}_i} \vert t_{j,1}x_1^{(i)} + t_{j,2}x_2^{(i)} + \ldots + t_{j,\ell} x^{(i)}_{\ell} \vert }{\min\limits_{j \in \mathcal{F}_i} |t_{j,1}x_1^{(i)} + t_{j,2}x_2^{(i)} + \ldots + t_{j,\ell} x^{(i)}_{\ell} |}.
\end{align}

For arbitrary $i = 1,\ldots,b$ with $\mT \in \mathbb{R}^{b \times \ell}$ and $\ell > 3^d+1$,  the null space of submatrix $\mT^{(i)}$ has dimension at least $s = \ell - 3^d > 1$. 
We compute
\begin{align*}
    \begin{bmatrix} \xx_1^{(i)} & \ldots & \xx_s^{(i)} \end{bmatrix} = \nullsp \left (\mT^{(i)},s \right ),
\end{align*}
and write any normalized vector in the null space as
\begin{align*}
    \xx = \alpha_1(\mtx{\theta}) \xx_1^{(i)} + \ldots + \alpha_s(\mtx{\theta}) \xx_s^{(i)}
\end{align*}
for (spherical) coordinates $\left (\alpha_1(\mtx{\theta}),\ldots,\alpha_s(\mtx{\theta}) \right ) \in \mathbb{R}^{s}$ parameterizing the unit hypersphere over $\mtx{\theta} \in \mathbb{R}^{s-1}$.
We then solve the small constrained optimization problem
\begin{align*}
    \mtx{\theta}^* = \arg \min_{\mtx{\theta} \in \mathbb{R}^{s-1}} \frac{ \max\limits_{j \in \mathcal{F}_i} \left | \left ( \alpha_1(\mtx{\theta}) \xx_1^{(i)} + \ldots + \alpha_s(\mtx{\theta}) \xx_s^{(i)} \right )^T \vct{t}^{(j)}   \right |  }{ \min\limits_{j \in \mathcal{F}_i} \left | \left ( \alpha_1(\mtx{\theta}) \xx_1^{(i)} + \ldots + \alpha_s(\mtx{\theta}) \xx_s^{(i)} \right )^T \vct{t}^{(j)}   \right |},
\end{align*}
where $\mtx{t}^{(j)} = \begin{bmatrix} t_{j,1} & t_{j,2} & \ldots & t_{j,\ell} \end{bmatrix}^T$,
to write our desired null space vector as 
\begin{align*}
    \vz^{(i)} =  \alpha_1(\mtx{\theta}^*) \xx_1^{(i)} + \ldots + \alpha_s(\mtx{\theta}^*) \xx_s^{(i)},
\end{align*}
and compute the basis matrices $\mU_{i}$, $\mV_i$ as in Section~\ref{sec:tagging}.

The overall computational cost of this optimization procedure is negligible when it is integrated into Algorithm~\ref{alg:stepI_tagging}, since the optimization happens over a convex region, with no more than a 2-dimensional parameterization in practice.
However, each additional column in the tagging matrix corresponds to $k+p$ additional samples of $\ma$ and $\ma^*$ with $\momega$ and $\mpsi$, though the total is still far fewer than in block nullification which we verify numerically in Section~\ref{sec:num_exp}. 
First, for completeness, we outline our method of reconstructing the full uniform BLR representation of (\ref{eq:uniformBLRfull}). 

\section{Randomized Compression Algorithms for Uniform BLR Matrices}
\label{sec:full_comp_algs}

In the previous sections, we focused on step~\ref{comp1} of randomized compression---the computation of basis matrices $\mU, \mV$---due to the similarity of compression algorithms after step~\ref{comp1}.
More precisely, steps~\ref{comp2} and \ref{comp3} of uniform BLR compression can be executed in a manner that is oblivious to the particular method used for step~\ref{comp1}, allowing for direct performance comparisons of block nullification (Algorithm~\ref{alg:stepI_bn}) and tagging (Algorithm \ref{alg:stepI_tagging}).
\textcolor{black}{For the sake of direct comparisons, we choose the same fixed target ranks for all basis computation methods in step~\ref{comp1}; however, in practice, more sophisticated techniques could be used instead.
For example, target ranks could be chosen adaptively using the adaptive randomized rangefinder, cf. Section~\ref{sec:preliminaries_randomized}, or larger target ranks could be chosen and truncated as needed.}

In this section, we briefly describe how steps~\ref{comp1}-\ref{comp3} are conducted in our experiments.
First, Algorithm~\ref{alg:stepI_basic} introduces the last basis construction algorithm we use as a benchmark for step~\ref{comp1}, which is equivalent to a blocked version of the basic randomized SVD with $\mathcal{O}(bk)$ structured Gaussian test matrices. 
We then describe a uniform BLR compression procedure that can be performed with basis matrices obtained from any of our algorithms; we discuss compression algorithms that reuse the sketches from step~\ref{comp1} for steps~\ref{comp2} and \ref{comp3} in the Supplementary Materials.


\begin{algorithm}[h]
    \caption{Basic RandSVD for Basis Construction}\label{alg:stepI_basic} \textbf{Call:} $[\mU,\mV] = \texttt{basic}(\ma, \ma^*, k)$
    \begin{algorithmic}[1]
        \Require Fast matrix-vector multiplication with uniform BLR $\ma \in \mathbb{R}^{N \times N}$ and $\ma^* \in \mathbb{R}^{N \times N}$, target rank $k$, given $b \times b$ tessellation of $d$-dimensional geometry with max block-size $m$
        \Ensure $\mU, \mV \in \mathbb{R}^{N \times bk}$ in uniform BLR representation of $\ma$ as in (\ref{eq:uniformBLRfull})

        \State Set $r = k+p$
        \For{blocks $i =1,\ldots,b$}
            \State Draw independent Gaussian test matrices $\momega, \mpsi \in \mathbb{R}^{N \times r}$
            \State Set $\momega(I_j,:) = 0$ and $\mpsi(I_j,:) = 0$ for all $j \in \mathcal{N}_i$
            \State Form $\mY = \ma \momega$ and $\mZ = \ma^* \mpsi$
            \State Compute $\mU_i = \col \left (\mY(I_i,:), \ k \right )$
            \State Compute $\mV_{i} = \col \left ( \mZ(I_i,:), \ k \right )$
        \EndFor
        \State Set $\mU = \textrm{diag} \left (\mU_1,\ldots,\mU_b \right )$ and $\mV = \textrm{diag} \left (\mV_1,\ldots,\mV_b \right ) $
    \end{algorithmic}
\end{algorithm}

\subsection{Direct evaluation in steps~\ref{comp2} and \ref{comp3}}
\label{sec:typeAalgs}

We present algorithms to compress uniform BLR matrices that allow for the most direct comparison of our basis construction algorithms.
First, we summarize in Algorithm~\ref{alg:stepI_basic} an algorithm used only as a benchmark for step~\ref{comp1}, which is equivalent to a blocked version of the randomized SVD done ``basically''  with $\mathcal{O}(bk)$ structured Gaussian test matrices.

Steps~\ref{comp2} and \ref{comp3} for (\ref{eq:uniformBLRfull}) can be performed straightforwardly by evaluating $\widetilde{\ma} = \mU^* (\ma \mV)$ directly, for $bk$ additional matvecs with $\ma$.
Then $3^d$ sparse structured matrices of size $N \times m$, containing the $m \times m$ identity as submatrices, can be used to extract the nonzero entries of $\ma - \mU \widetilde{\ma} \mV^*$ for $\mB$ as in \cite{lin2011fast}, e.g.

\begin{center}
\begin{minipage}{0.7\linewidth}
\begin{align}
\begin{array}{cccc}
\mY_{\mB} & = & (\ma - \mU \widetilde{\ma} \mV^*) & \momega_{\mB} \\
 \begin{tikzpicture}[scale=0.63, line width=1pt]
\def\n{7cm}
            \draw[fill=\lrcolor, dashed] (\n*0/8, \n*1/8) rectangle (\n*1/8, \n*6/8) node[midway] {};
            \draw[fill=\lrcolor, dashed] (\n*0/8, \n*2/8) rectangle (\n*1/8, \n*3/8) node[midway] {};
            \draw[fill=\lrcolor] (\n*0/8, \n*0/8) rectangle (\n*1/8, \n*1/8) node[midway] {\scriptsize $\mB_{7,8}$};          
            \draw[fill=\lrcolor] (\n*0/8, \n*7/8) rectangle (\n*1/8, \n*8/8) node[midway] {\scriptsize  $\mB_{1,1}$};
            \draw[fill=\lrcolor] (\n*0/8, \n*4/8) rectangle (\n*1/8, \n*5/8) node[midway] {\scriptsize  $\mB_{4,4}$};
            \draw[fill=\lrcolor] (\n*0/8, \n*1/8) rectangle (\n*1/8, \n*2/8) node[midway] {\scriptsize $\mB_{7,7}$};
            \draw[fill=\lrcolor] (\n*0/8, \n*6/8) rectangle (\n*1/8, \n*7/8) node[midway] {\scriptsize $\mB_{2,1}$};
            \draw[fill=\lrcolor] (\n*0/8, \n*3/8) rectangle (\n*1/8, \n*4/8) node[midway] {\scriptsize $\mB_{5,4}$};
\end{tikzpicture} &  & \fontsize{8}{\baselineskip}
    \def\n{7cm}
    \begin{tikzpicture}[level distance=\n/10,
        level 1/.style={sibling distance=\n/2},
        level 2/.style={sibling distance=\n/4},
        level 3/.style={sibling distance=\n/8},
        scale=0.63,]
        \begin{scope}[local bounding box=scope1, line width=1pt]

            \draw (0, 0) rectangle (\n, \n);

            \draw[fill=white] (\n*2/8, \n*7/8) rectangle (\n*3/8, \n*8/8) node[midway] {};
            \draw[fill=white] (\n*3/8, \n*7/8) rectangle (\n*4/8, \n*8/8) node[midway] {};
            \draw[fill=white] (\n*4/8, \n*7/8) rectangle (\n*5/8, \n*8/8) node[midway] {};
            \draw[fill=white] (\n*5/8, \n*7/8) rectangle (\n*6/8, \n*8/8) node[midway] {};
            \draw[fill=white] (\n*6/8, \n*7/8) rectangle (\n*7/8, \n*8/8) node[midway] {};
            \draw[fill=white] (\n*7/8, \n*7/8) rectangle (\n*8/8, \n*8/8) node[midway] {};

            \draw[fill=white] (\n*3/8, \n*6/8) rectangle (\n*4/8, \n*7/8) node[midway] {};
            \draw[fill=white] (\n*4/8, \n*6/8) rectangle (\n*5/8, \n*7/8) node[midway] {};
            \draw[fill=white] (\n*5/8, \n*6/8) rectangle (\n*6/8, \n*7/8) node[midway] {};
            \draw[fill=white] (\n*6/8, \n*6/8) rectangle (\n*7/8, \n*7/8) node[midway] {};
            \draw[fill=white] (\n*7/8, \n*6/8) rectangle (\n*8/8, \n*7/8) node[midway] {};

            \draw[fill=white] (\n*0/8, \n*5/8) rectangle (\n*1/8, \n*6/8) node[midway] {};
            \draw[fill=white] (\n*4/8, \n*5/8) rectangle (\n*5/8, \n*6/8) node[midway] {};
            \draw[fill=white] (\n*5/8, \n*5/8) rectangle (\n*6/8, \n*6/8) node[midway] {};
            \draw[fill=white] (\n*6/8, \n*5/8) rectangle (\n*7/8, \n*6/8) node[midway] {};
            \draw[fill=white] (\n*7/8, \n*5/8) rectangle (\n*8/8, \n*6/8) node[midway] {};

            \draw[fill=white] (\n*0/8, \n*4/8) rectangle (\n*1/8, \n*5/8) node[midway] {};
            \draw[fill=white] (\n*1/8, \n*4/8) rectangle (\n*2/8, \n*5/8) node[midway] {};
            \draw[fill=white] (\n*5/8, \n*4/8) rectangle (\n*6/8, \n*5/8) node[midway] {};
            \draw[fill=white] (\n*6/8, \n*4/8) rectangle (\n*7/8, \n*5/8) node[midway] {};
            \draw[fill=white] (\n*7/8, \n*4/8) rectangle (\n*8/8, \n*5/8) node[midway] {};

            \draw[fill=white] (\n*0/8, \n*3/8) rectangle (\n*1/8, \n*4/8) node[midway] {};
            \draw[fill=white] (\n*1/8, \n*3/8) rectangle (\n*2/8, \n*4/8) node[midway] {};
            \draw[fill=white] (\n*2/8, \n*3/8) rectangle (\n*3/8, \n*4/8) node[midway] {};
            \draw[fill=white] (\n*6/8, \n*3/8) rectangle (\n*7/8, \n*4/8) node[midway] {};
            \draw[fill=white] (\n*7/8, \n*3/8) rectangle (\n*8/8, \n*4/8) node[midway] {};

            \draw[fill=white] (\n*0/8, \n*2/8) rectangle (\n*1/8, \n*3/8) node[midway] {};
            \draw[fill=white] (\n*1/8, \n*2/8) rectangle (\n*2/8, \n*3/8) node[midway] {};
            \draw[fill=white] (\n*2/8, \n*2/8) rectangle (\n*3/8, \n*3/8) node[midway] {};
            \draw[fill=white] (\n*3/8, \n*2/8) rectangle (\n*4/8, \n*3/8) node[midway] {};
            \draw[fill=white] (\n*7/8, \n*2/8) rectangle (\n*8/8, \n*3/8) node[midway] {};

            \draw[fill=white] (\n*0/8, \n*1/8) rectangle (\n*1/8, \n*2/8) node[midway] {};
            \draw[fill=white] (\n*1/8, \n*1/8) rectangle (\n*2/8, \n*2/8) node[midway] {};
            \draw[fill=white] (\n*2/8, \n*1/8) rectangle (\n*3/8, \n*2/8) node[midway] {};
            \draw[fill=white] (\n*3/8, \n*1/8) rectangle (\n*4/8, \n*2/8) node[midway] {};
            \draw[fill=white] (\n*4/8, \n*1/8) rectangle (\n*5/8, \n*2/8) node[midway] {};

            \draw[fill=white] (\n*0/8, \n*0/8) rectangle (\n*1/8, \n*1/8) node[midway] {};
            \draw[fill=white] (\n*1/8, \n*0/8) rectangle (\n*2/8, \n*1/8) node[midway] {};
            \draw[fill=white] (\n*2/8, \n*0/8) rectangle (\n*3/8, \n*1/8) node[midway] {};
            \draw[fill=white] (\n*3/8, \n*0/8) rectangle (\n*4/8, \n*1/8) node[midway] {};
            \draw[fill=white] (\n*4/8, \n*0/8) rectangle (\n*5/8, \n*1/8) node[midway] {};
            \draw[fill=white] (\n*5/8, \n*0/8) rectangle (\n*6/8, \n*1/8) node[midway] {};

            \draw[fill=\frcolor] (\n*1/8, \n*7/8) rectangle (\n*2/8, \n*8/8) node[midway] {};
            \draw[fill=\frcolor] (\n*2/8, \n*6/8) rectangle (\n*3/8, \n*7/8) node[midway] {};
            \draw[fill=\frcolor] (\n*3/8, \n*5/8) rectangle (\n*4/8, \n*6/8) node[midway] {};
            \draw[fill=\frcolor] (\n*4/8, \n*4/8) rectangle (\n*5/8, \n*5/8) node[midway] {};
            \draw[fill=\frcolor] (\n*5/8, \n*3/8) rectangle (\n*6/8, \n*4/8) node[midway] {};
            \draw[fill=\frcolor] (\n*6/8, \n*2/8) rectangle (\n*7/8, \n*3/8) node[midway] {};
            \draw[fill=\frcolor] (\n*7/8, \n*1/8) rectangle (\n*8/8, \n*2/8) node[midway] {};

            \draw[fill=\frcolor] (\n*0/8, \n*6/8) rectangle (\n*1/8, \n*7/8) node[midway] {};
            \draw[fill=\frcolor] (\n*1/8, \n*5/8) rectangle (\n*2/8, \n*6/8) node[midway] {};
            \draw[fill=\frcolor] (\n*2/8, \n*4/8) rectangle (\n*3/8, \n*5/8) node[midway] {};
            \draw[fill=\frcolor] (\n*3/8, \n*3/8) rectangle (\n*4/8, \n*4/8) node[midway] {};
            \draw[fill=\frcolor] (\n*4/8, \n*2/8) rectangle (\n*5/8, \n*3/8) node[midway] {};
            \draw[fill=\frcolor] (\n*5/8, \n*1/8) rectangle (\n*6/8, \n*2/8) node[midway] {};
            \draw[fill=\frcolor] (\n*6/8, \n*0/8) rectangle (\n*7/8, \n*1/8) node[midway] {};

            \draw[fill=\frcolor] (\n*0/8, \n*7/8) rectangle (\n*1/8, \n*8/8) node[midway] {};
            \draw[fill=\frcolor] (\n*1/8, \n*6/8) rectangle (\n*2/8, \n*7/8) node[midway] {};
            \draw[fill=\frcolor] (\n*2/8, \n*5/8) rectangle (\n*3/8, \n*6/8) node[midway] {};
            \draw[fill=\frcolor] (\n*3/8, \n*4/8) rectangle (\n*4/8, \n*5/8) node[midway] {};
            \draw[fill=\frcolor] (\n*4/8, \n*3/8) rectangle (\n*5/8, \n*4/8) node[midway] {};
            \draw[fill=\frcolor] (\n*5/8, \n*2/8) rectangle (\n*6/8, \n*3/8) node[midway] {};
            \draw[fill=\frcolor] (\n*6/8, \n*1/8) rectangle (\n*7/8, \n*2/8) node[midway] {};
            \draw[fill=\frcolor] (\n*7/8, \n*0/8) rectangle (\n*8/8, \n*1/8) node[midway] {};
        \end{scope}
    \end{tikzpicture} & \begin{tikzpicture}[scale=1.0, line width=1pt]
            \def\n{4cm}
            \draw[fill=\frcolor] (\n*0/8, \n*7/8) rectangle (\n*2/8, \n*8/8) node[midway] {$\mtx{I}$};
            \draw[fill=white] (\n*0/8, \n*6/8) rectangle (\n*2/8, \n*7/8) node[midway] {};
            \draw[fill=white] (\n*0/8, \n*5/8) rectangle (\n*2/8, \n*6/8) node[midway] {};
            \draw[fill=\frcolor] (\n*0/8, \n*4/8) rectangle (\n*2/8, \n*5/8) node[midway] {$\mtx{I}$};
            \draw[fill=white] (\n*0/8, \n*3/8) rectangle (\n*2/8, \n*4/8) node[midway] {};
            \draw[fill=white] (\n*0/8, \n*2/8) rectangle (\n*2/8, \n*3/8) node[midway] {};
            \draw[fill=\frcolor] (\n*0/8, \n*1/8) rectangle (\n*2/8, \n*2/8) node[midway] {$\mtx{I}$};
            \draw[fill=white] (\n*0/8, \n*0/8) rectangle (\n*2/8, \n*1/8) node[midway] {};
\end{tikzpicture} 
\end{array}  .  
\label{eq:sparse_struct_ids}
\end{align}
\end{minipage}
\end{center}


\textcolor{black}{Table~\ref{tab:full_comp_algs_A} summarizes each reconstruction algorithm we consider using the basis functions \texttt{bn} (Algorithm~\ref{alg:stepI_bn}), \texttt{tag} (Algorithm~\ref{alg:stepI_tagging}), and \texttt{basic} (Algorithm~\ref{alg:stepI_basic}).
Steps 2--3 in Table~\ref{tab:full_comp_algs_A} are all identical, making this reconstruction procedure ideal for performance comparisons of our basis computation algorithms.}

\begin{table}[h]
\centering
\begin{tabular}{lll}
  {\texttt{reconstruct\_bn}} & {\texttt{reconstruct\_tag}} & {\texttt{reconstruct\_basic}} \\
  \begin{minipage}{0.28\linewidth}\begin{enumerate}[wide, labelwidth=!,itemindent=!,labelindent=0pt, leftmargin=0em, label=\arabic*., itemsep=1mm, parsep=0pt]
  \vspace{1.2mm}
    \item $[\mU,\mV] = \texttt{bn}(\ma,\ma^*,k)$
    \item Form $\widetilde{\ma} = \mU^*(\ma \mV)$ 
    \item Form $\mB $ as in (\ref{eq:sparse_struct_ids})
  \end{enumerate}\end{minipage} &
  \begin{minipage}{0.3\linewidth}\begin{enumerate}[wide, labelwidth=!,itemindent=!,labelindent=0pt, leftmargin=0em, label=\arabic*., itemsep=1mm, parsep=0pt]
  \vspace{1.2mm}
    \item $[\mU,\mV] = \texttt{tag}(\ma,\ma^*,k)$
    \item Form $\widetilde{\ma} = \mU^*(\ma \mV)$ 
    \item Form $\mB $ as in (\ref{eq:sparse_struct_ids})
  \end{enumerate}\end{minipage} &
  \begin{minipage}{0.33\linewidth}\begin{enumerate}[wide, labelwidth=!,itemindent=!,labelindent=0pt, leftmargin=0em, label=\arabic*., itemsep=1mm, parsep=0pt]
  \vspace{1.2mm}
    \item $[\mU,\mV] = \texttt{basic}(\ma,\ma^*,k)$
    \item Form $\widetilde{\ma} = \mU^*(\ma \mV)$ 
    \item Form $\mB $ as in (\ref{eq:sparse_struct_ids})
  \end{enumerate}\end{minipage}
\end{tabular}
 \caption{\textcolor{black}{Randomized compression algorithms for uniform BLR format of (\ref{eq:uniformBLRfull}). Steps 2--3 contribute an additional $bk$ and $3^dm$ matvecs with $\ma$, respectively. }}
\label{tab:full_comp_algs_A}
\end{table}

\section{Numerical Experiments}
\label{sec:num_exp}

We now demonstrate the improved performance of tagging over block nullification to compress strongly admissible uniform BLR matrices.
For several different test problems and problem sizes $N$, we report the following:
\begin{itemize}
    \item Accuracy of compressed $\ma_{\textup{uBLR}}$ of the form (\ref{eq:uniformBLRfull}), using the relative error metric $\frac{\|\ma - \ma_{\textup{uBLR}} \|_2}{\|\ma\|_2}$ via 20 iterations of the randomized power method \cite{halko2011finding},
    \item Total runtime (in seconds) of each compression algorithm, 
    \item Total number of matvecs with $\ma$, $\ma^*$ required for compression. 
\end{itemize}
For tagging, we also report the aspect ratios after performing the optimization routine of Section~\ref{sec:tagmat_practical}. \textcolor{black}{All experiments were performed with Gaussian tagging matrices with 3 additional columns; alternative tagging matrices are discussed in the Supplementary Materials. More aggressive numbers of extra columns were tested experimentally, but provided little benefit after 3.}
All test problems were implemented in MATLAB 2024a, and all experiments were carried out on a workstation with an Intel(R) Xeon(R) Gold 6254 CPU operating at 3.10GHz with 72 cores and 750 GB of memory.

For each test problem, we use a target block-rank of $k=30$ with an oversampling parameter of $p=10$. 
We also report the number of blocks $b$ and the maximum block size $m$ for each problem size $N$.
To substantiate our choices of $b$ and $m$, we highlight a key distinction between hierarchical and flat rank-structured matrix formats. 
Often in hierarchical rank-structured matrix compression, the leaf node size $m$ is chosen such that $m = O(k+p)$ with approximately $2N/m$ total nodes in the index tree; moreover, linear complexity can be achieved by leveraging nested bases, e.g. \cite{levitt2022linear}.

By contrast, randomized compression of flat formats does not attain linear complexity.
The storage requirement of a strongly admissible uniform BLR matrix is 
$M\ \sim\ Nk + b^2 k^2 + 3^d {N^2}/b$ floating point values,
and the compressed uniform BLR representation $\ma_{\textup{uBLR}}$ of (\ref{eq:uniformBLRfull})
can be recovered in no fewer than
$N_{\rm matvec} \sim M / N$ matvecs, 
since a matrix of size $N \times N_{\rm matvec}$ holds the \textit{minimum} number of entries needed to
store $\ma_{\textup{uBLR}}$. 
The dominant storage costs can be attributed to $\widetilde{\ma}$ and $\mB$, and it is challenging to recover them in an ``optimal'' number of matvecs ($N_{\rm matvec} \sim 3^d m + bk$).
Thus, we choose $b$ to balance $N_{\rm matvec}$ so that $b = \sqrt{3^dN/k} \ \Rightarrow \ M \sim \sqrt{ 3^dk }\ N^{3/2},$
a reasonable choice for medium-sized problems, e.g. $20,000 \leq N \leq 100,000$ here.
\textcolor{black}{In other words, in our experiments, $b$ ranges from $77$ to $173$ for $N = 20,000$ to $100,000$, with $m = N/b$}.

\subsection{2D Laplace Kernel}\label{sub:2d_laplace_kernel}

To profile performance, we use the Green's function of 
the 2D Laplace equation for uniformly random $\{x_i\}_{i=1}^N$ in the unit square, where 
\begin{equation}
\mtx A_{ij} = \log(\|x_i - x_j\|),\ \text{for}\ i \neq j,
\end{equation}
and entries on the diagonal are set to 0. Dense systems of this form commonly arise in the context of integral equations.
In practice, the method of proxy surfaces is often used to approximate basis matrices algebraically \cite{cheng2005compression,xing2020interpolative}, but we include this problem as a benchmark because the algebraic rank behavior of $\mtx A$ is well-characterized by multipole estimates \cite{greengard1987fast,greengard1997new} and exhibits exponential decay. 

\begin{figure}[htb!]

\begin{subfigure}{0.33\linewidth}
    \centering
    \includegraphics[width=\linewidth]{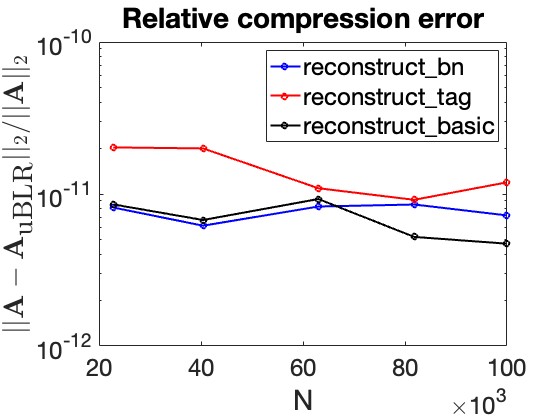}
    \caption{Relative error.}
    \label{fig:fmm_rel_error}
\end{subfigure}
\begin{subfigure}{0.33\linewidth}
    \centering
    \includegraphics[width=\linewidth]{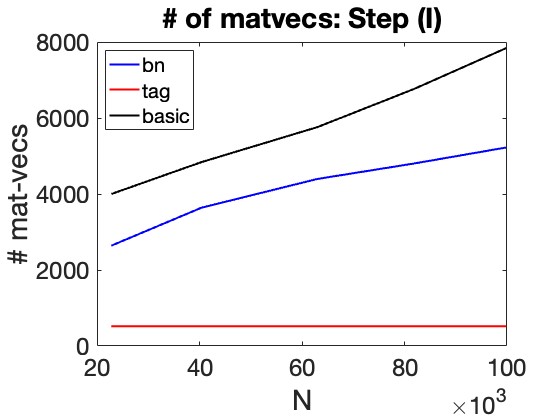}
    \caption{Step~\ref{comp1} matvecs}
    \label{fig:fmm_basis_mv}
\end{subfigure}
\begin{subfigure}{0.32\linewidth}
    \centering
    \includegraphics[width=\linewidth]{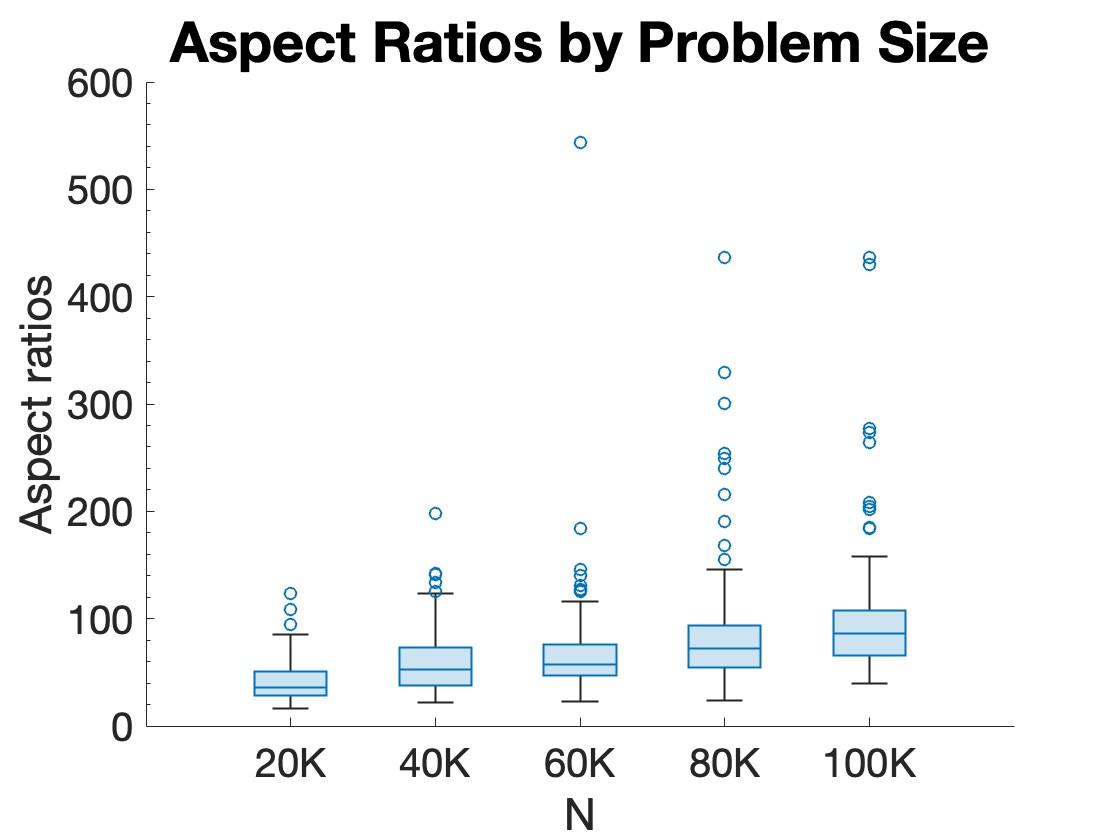}
    \caption{Tagging aspect ratios.}
    \label{fig:fmm_tagging_ar}
\end{subfigure}

\vskip 0.25em

\begin{subfigure}{0.33\linewidth}
    \centering
    \includegraphics[width=\linewidth]{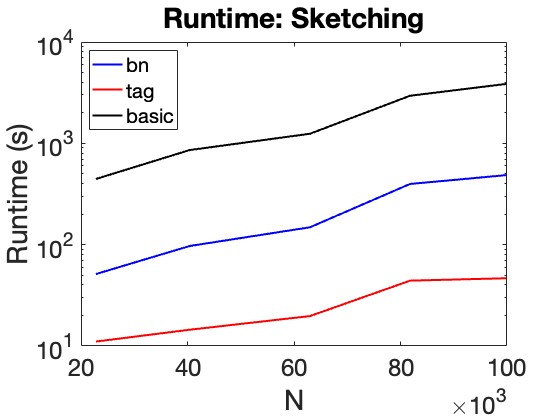}
    \caption{Sketching runtimes.}
    \label{fig:fmm_sketching_time}
\end{subfigure}%
\begin{subfigure}{0.33\linewidth}
    \centering
    \includegraphics[width=\linewidth]{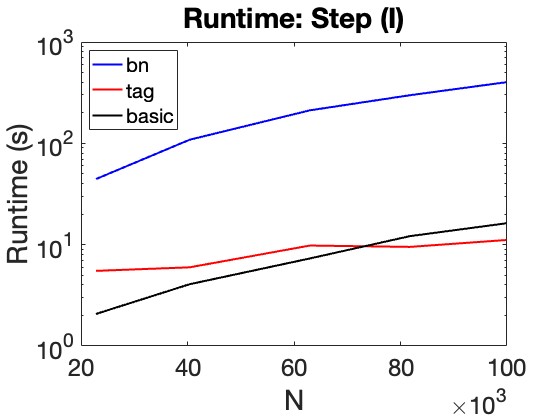}
    \caption{Step~\ref{comp1} runtimes.}
    \label{fig:fmm_basis_time}
\end{subfigure}
\begin{subfigure}{0.33\linewidth}
    \centering
    \includegraphics[width=\linewidth]{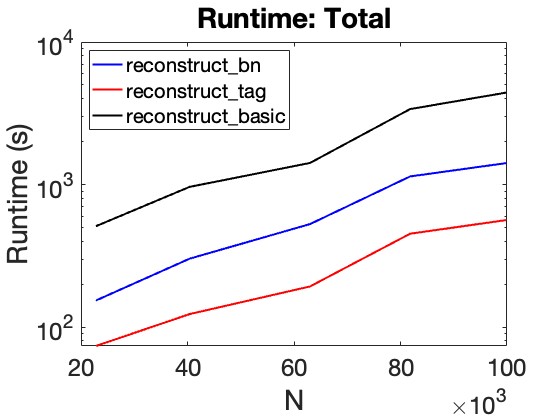}
    \caption{Total runtimes.}
    \label{fig:fmm_total_time}
\end{subfigure}

\caption{Accuracy, number of matvecs, and timing results for 2D Laplace FMM in Section~\ref{sub:2d_laplace_kernel}. 
Fig. \ref{fig:fmm_rel_error} reports the reconstruction accuracy for algorithmic variants.
Fig.~\ref{fig:fmm_basis_mv} reports the number of matrix-vector products for step~\ref{comp1}.
Fig.~\ref{fig:fmm_tagging_ar} shows box-plots of the tagging aspect ratios.
Figs.~\ref{fig:fmm_sketching_time}-\ref{fig:fmm_total_time} report algorithmic runtimes.}
\label{fig:fmm_results}
\end{figure}

Fig.~\ref{fig:fmm_results} summarizes the results. Matrix-vector products \(\mtx{x} \rightarrow \mtx{A} \mtx{x}\) are performed using \texttt{FMM2D}, a Fortran implementation of the fast multipole method developed and maintained by the Flatiron Institute. Key observations are as follows:

\begin{itemize}
    \item Tagging significantly reduces the number of matrix-vector products (\(N_{\rm matvec}\)) required for basis construction compared to alternative methods. \textcolor{black}{For \texttt{bn} (Algorithm~\ref{alg:stepI_bn}) and \texttt{basic} (Algorithm~\ref{alg:stepI_basic}), \(N_{\rm matvec}\) scales with the block size and the number of blocks, respectively--both of which grow with the problem size \(N\) in flat formats. In contrast, the number of matvecs for \texttt{tag} (Algorithm~\ref{alg:stepI_tagging}) remains constant.}

    \item The reduction in \(N_{\rm matvec}\) for basis construction translates into considerable time savings during the sketching of \(\mtx{Y}\) and \(\mtx{Z}\).

    \item Tagging saves substantially on post-processing time vs. block nullification.

    \item \textcolor{black}{The full algorithm \texttt{reconstruct\_tag} achieves an 8.3x reduction in the total number of matvecs for the largest problem size compared to naive matrix formation.}

    \item Tagging maintains comparable accuracy to alternative methods, and the aspect ratios can be effectively controlled by the method of Section~\ref{sec:tagmat_practical}.
\end{itemize}

\subsection{Sparse Schur Complement for a Thin Slab}
\label{sec:schur}
In sparse direct solvers for elliptic PDEs, compressing and factorizing sparse matrices is often necessary. Accessing matrix entries directly is computationally challenging, and randomized sketching techniques are frequently used to accelerate and simplify nested dissection solvers.

Consider solving the constant-coefficient Helmholtz equation with zero body load and prescribed Dirichlet boundary conditions on a domain $\Omega$:
\begin{equation}
\label{eq:bvp}
\begin{aligned}
- \Delta u(x) - \kappa^2 u(x) &= 0, \quad & x \in \Omega, \\
u(x) &= g(x), \quad & x \in \partial \Omega.
\end{aligned}
\end{equation}
Discretizing with second-order finite differences leads to the linear system 
$\mtx{A} \mtx{u} = \mtx{f}$ to solve.
To solve this system efficiently, the domain $\Omega$ is partitioned into \textit{thin slabs}, where one dimension is constrained to be electrically small, as demonstrated in \cite{yesypenko2024slablu,engquist2011sweeping}.

\begin{figure}[!htb]
\centering
\begin{minipage}{0.6\textwidth}
    \centering    \includegraphics[width=0.6\textwidth]{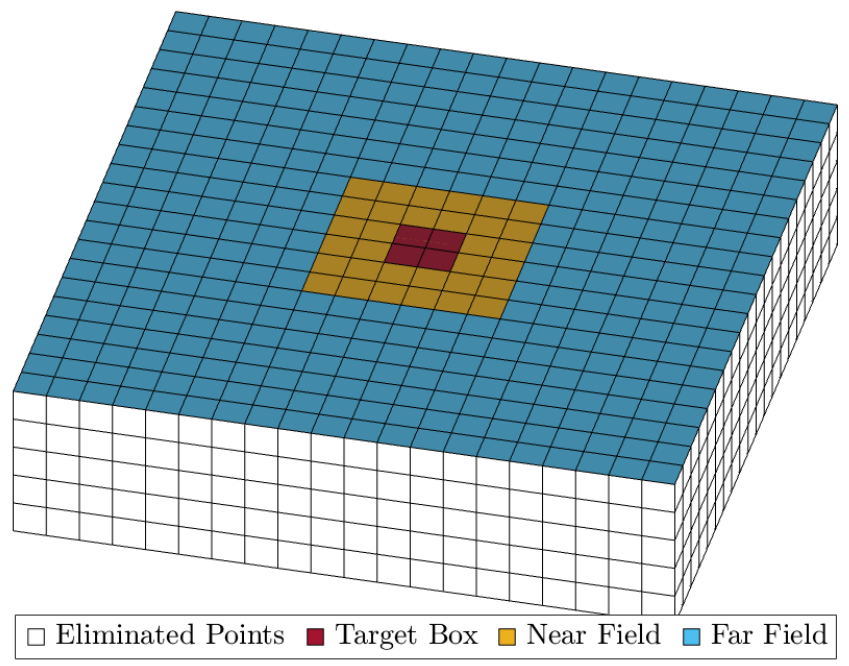}
\end{minipage}%
\hfill
\begin{minipage}{0.4\textwidth}
    \caption{The discretization of a thin slab of size $n \times n \times l$, where $l \ll n$ is fixed to be $l=10$ in our experiments. The eliminated discretization points are shown in white, and the frontal indices are colored to denote the far field and near field for a target box.}
    \label{fig:slab_discretization}
\end{minipage}
\end{figure}

We investigate the use of uniform BLR matrices for a slab subdomain. For a domain with $n \times n \times l$ discretization points, where $l$ is fixed to be $l=10$, the front size grows as $N = n^2$ as visualized in Fig.~\ref{fig:slab_discretization}. The wavenumber parameter scales with the number of discretization points to maintain 100 points per wavelength. For the largest problem size ($N = 100,000$), the domain measures approximately $3\lambda \times 3\lambda \times 0.1\lambda$, where $\lambda$ denotes the wavelength.
In the context of domain decomposition, the degrees of freedom are partitioned into frontal nodes and internal nodes. 
The Schur complement is the linear algebraic operator that eliminates the internal nodes in a multifrontal solver, resulting in a dense matrix defined on the frontal nodes.
While the Schur complement is dense, it can be applied efficiently to vectors by leveraging the sparsity of its components \cite{yesypenko2024slablu}, enabling fast application of the solver to vectors.

\begin{figure}[htb!]

\begin{subfigure}{0.33\linewidth}
    \centering
    \includegraphics[width=\linewidth]{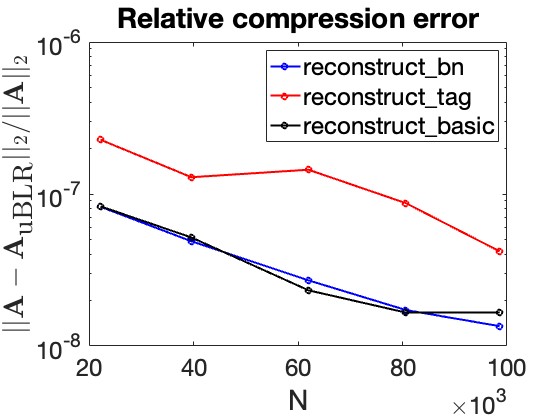}
    \caption{Relative error.}
    \label{fig:schur_rel_error}
\end{subfigure}
\begin{subfigure}{0.33\linewidth}
    \centering
    \includegraphics[width=\linewidth]{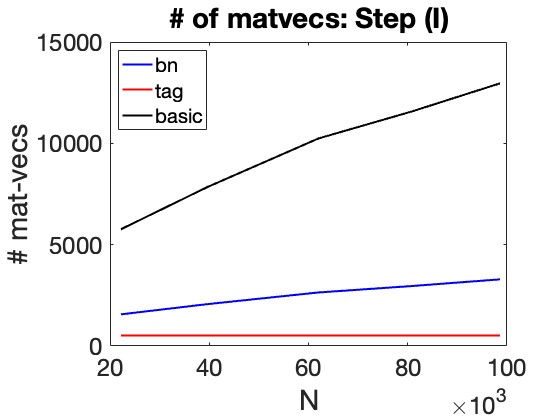}
    \caption{Step~\ref{comp1} matvecs.}
    \label{fig:schur_basis_mv}
\end{subfigure}
\begin{subfigure}{0.32\linewidth}
    \centering
    \includegraphics[width=\linewidth]{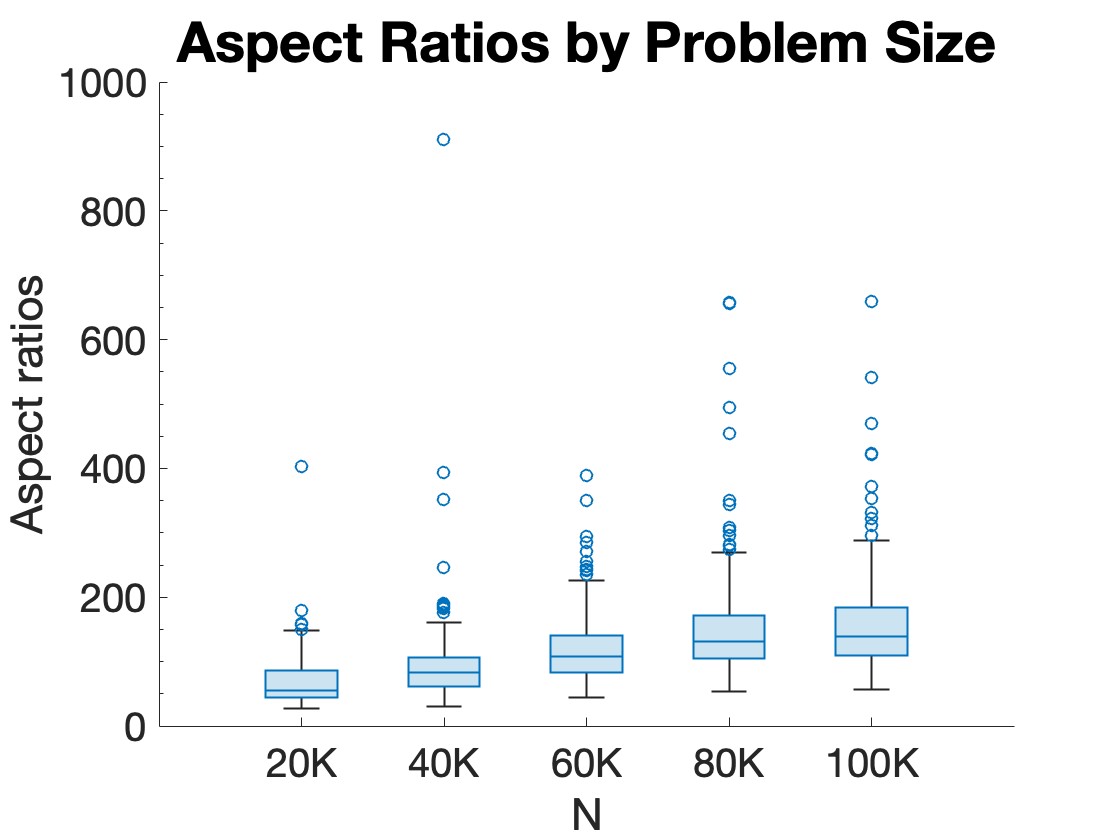}
    \caption{Tagging aspect ratios.}
    \label{fig:schur_tag_ar}
\end{subfigure}

\vskip 0.25em

\begin{subfigure}{0.33\linewidth}
    \centering
    \includegraphics[width=\linewidth]{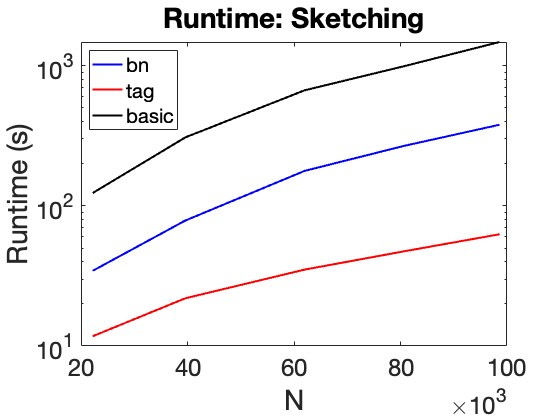}
    \caption{Sketching runtimes.}
    \label{fig:schur_sketching_time}
\end{subfigure}%
\begin{subfigure}{0.33\linewidth}
    \centering
    \includegraphics[width=\linewidth]{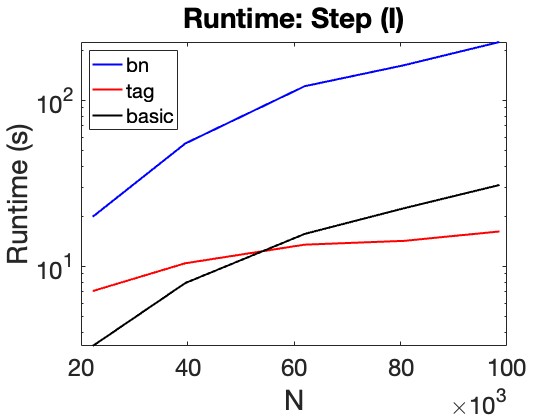}
    \caption{Step~\ref{comp1} runtimes.}
    \label{fig:schur_basis_time}
\end{subfigure}
\begin{subfigure}{0.33\linewidth}
    \centering
    \includegraphics[width=\linewidth]{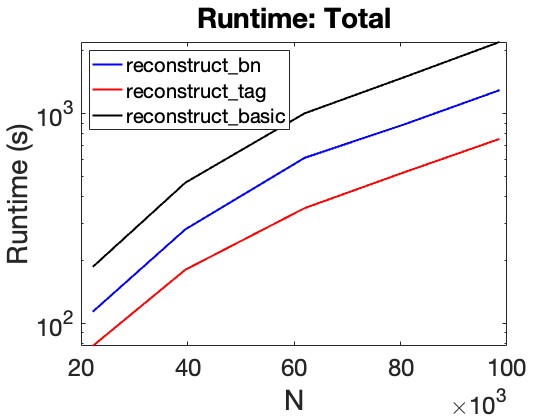}
    \caption{Total runtimes.}
    \label{fig:schur_total_time}
\end{subfigure}

\caption{Accuracy, number of matvecs, and timing results for the Schur complement in Section~\ref{sec:schur}. 
Fig. \ref{fig:schur_rel_error} reports the reconstruction accuracy for algorithmic variants.
Fig.~\ref{fig:schur_basis_mv} reports the number of matrix-vector products for step~\ref{comp1}.
Fig.~\ref{fig:schur_tag_ar} shows box-plots of tagging aspect ratios.
Figs.~\ref{fig:schur_sketching_time}-\ref{fig:schur_total_time} report algorithmic runtimes.}
\label{fig:schur_results}
\end{figure}

Fig.~\ref{fig:schur_results} shows that tagging provides excellent scaling in the number of samples needed for step (I). It is also the most computationally efficient. As in Section \ref{sub:2d_laplace_kernel}, we observe that the far-field rank decays exponentially. Since the slab width is fixed, the rank of far-field interactions decays \textit{faster} as the problem size grows, leading to improved approximation accuracy for increasing $N$ and a fixed rank $k$.
As in Section~\ref{sub:2d_laplace_kernel}, the tagging aspect ratios can be effectively controlled by the method of Section~\ref{sec:tagmat_practical}.

\section{Conclusions and Future Work}
\label{sec:conclusions}

We present a black-box randomized compression algorithm based on tagging, improving on existing randomized compression algorithms for uniform BLR matrices under a strong admissibility condition.
Our method only requires $\mathcal{O}(k)$ random samples for basis computations, versus $\mathcal{O}(m+k)$ for block nullification (which increases with problem size $N$ for flat formats).
We show that compression with tagging achieves comparable accuracy to existing compression algorithms with greatly improved computational efficiency. 
We also draw a connection between optimality in tagging and Pl\"ucker coordinates, and we present an alternative fast method of optimizing tagging matrices that is reliable in practice. 

Avenues of future work include a hybrid numeric-symbolic computational scheme to generate theoretically optimal tagging matrix entries.
Additionally, we are developing a high-performance implementation of our randomized compression algorithm with tagging.
Future work will also investigate tagging for hierarchical formats.
{Thorough analysis of tagging is of interest as well, similar to the method presented in \cite{amsel2025quasioptimalhierarchicallysemiseparablematrix} which appeared on arXiv subsequent to this work.}

\section*{Acknowledgements}
The authors wish to thank Joe Kileel and Steven Karp for their valuable insights. 

\bibliographystyle{siamplain}
\bibliography{references}

\end{document}